\documentclass[12pt,a4paper]{article}

\usepackage{graphicx}
\usepackage{amssymb}
\usepackage{amsmath}
\usepackage{caption}
\usepackage{subcaption}
\usepackage[margin=0.8in]{geometry} % for the margins (for a 12pt document, the default margin is 1.5 in)
\usepackage{subfloat}
\usepackage{color}

\title{Analytical solutions for two-dimensional singly periodic Stokes flow singularity arrays near walls}

\author{Darren Crowdy$^{1}$ and Elena Luca$^{2}$}

\date{}

\usepackage{MnSymbol}
\usepackage{bigints}

\begin{document}

\maketitle

\begin{center}
$^{1}$Department of Mathematics \\
Imperial College London \\
180 Queen's Gate \\
London, SW7 2AZ, U.K.
\end{center}

\begin{center}
$^{2}$Department of Mechanical and Aerospace Engineering \\
Jacobs School of Engineering, UCSD \\
La Jolla, CA 92093-0411, USA.

\vskip 0.1truein 
{\tt d.crowdy@imperial.ac.uk} \\
{\tt elouca@eng.ucsd.edu}
\end{center}

\vskip 0.5truein

\begin{center}
{\bf Abstract}
\end{center}
\noindent
New analytical representations of the Stokes flows due to 
periodic arrays of point singularities in a two-dimensional no-slip channel and
 in the half-plane near a no-slip wall are derived. 
The analysis makes use of a conformal mapping from a concentric
annulus (or a disc) to a rectangle and a complex variable formulation of Stokes flow to derive the solutions. The form of the solutions is amenable to fast and accurate numerical computation without the need for Ewald summation or other fast summation techniques.

\vfill\eject

\section{Introduction}

There has been a resurgence of interest in the mathematical theory of Stokes flows as a result of the multifarious problems arising in burgeoning fields such as microfluidics, nonlinear electrokinetics, low-Reynolds-number swimming, superhydrophobic surfaces and the study of active suspensions. In many situations singularity theory is used to model complicated bodies, such as swimming microorganisms or colloidal particles: a low-Reynolds-number swimmer, which is force- and torque-free, is often modelled as a point stresslet \cite{Or}; Jeong \& Moffatt \cite{Jeong} modelled a pair of counter-rotating rollers beneath a free surface in a viscous fluid by an irrotational point dipole. Suspensions, or swarms of microorganisms, require the introduction of many singularities to model the complicated interaction mechanisms and it is common to consider flows comprising periodic arrays of the fundamental singularities of Stokes flow. 

Hasimoto \cite{Has2} was the first to consider doubly and triply periodic arrays of Stokes singularities,  recognising even then that the numerical computation of such flows is not without difficulty. He introduced the idea of using Ewald summation techniques to improve the numerical efficacy of computing such periodic flows. Pozrikidis \cite{Poz3} reconsidered similar problems from the perspective of constructing the basic Green's functions needed in formulating boundary integral methods for spatially periodic Stokes flows. In an earlier study, Pozrikidis \cite{Poz} used standard Fourier transform techniques to solve for a periodic array of point singularities in a two-dimensional channel and in the half-plane. Davis \cite{Davis1993} used Fourier transform techniques to solve various problems involving periodic arrays of Stokeslets in a two-dimensional channel and in the half-plane, and investigated the blocking properties of periodic arrays of wall-attached barriers. Work on improving numerical methods for Stokes flows in multi-particle settings continues \cite{Tornberg,Greengard} and extensions of these ideas to multi-particle interactions in confined geometries have also been made \cite{Hernandez}. Many of
these methods rely on a splitting of the flow into a local contribution for which rapidly decaying free-space analytical solutions can be employed together with a global contribution whose effect is determined numerically in an efficient way (either spectrally, or using an iterative scheme). 

The present authors \cite{CrowdyLuca2018} have recently given analytical representations 
that are amenable to fast numerical evaluation of the flows associated with doubly periodic arrangements of point singularities of two-dimensional Stokes flow.
They used analytic function theory, a conformal mapping from a concentric annulus
and the so-called Schottky--Klein prime function associated with that annulus to
derive their new form of the solutions.

The present paper makes a basic theoretical contribution to the study of {\em singly}
periodic two-dimensional Stokes flows in a no-slip channel and in the half-plane
near a single no-slip wall.
This paper can be viewed as a sequel to \cite{CrowdyLuca2018} where doubly periodic
arrays of Stokes flow singularities are considered.
Here we focus on singly periodic flows and, specifically, those
generated by arrays of fundamental singularities, i.e.~the Stokeslet and the higher order singularities. By employing an approach combining conformal mapping with a complex variable formulation of Stokes flow, we show how a fast and accurate 
representation of such flows can be derived without the need for the aforementioned ``splitting'' of the flow into separate local and global contributions or any Ewald summation techniques. While our results are limited to two dimensions, we believe they are valuable additions to the basic mathematical theory of Stokes flow.

\section{Stokes flows in two dimensions}

Consider a region of incompressible fluid  of viscosity $\eta$ governed by the
Stokes equations
\begin{equation}
\nabla p=\eta \nabla^{2} {\bf u}, \quad \nabla \cdot {\bf u}=0,
\end{equation}
where ${\bf u}=(u,v)$ is the two-dimensional velocity field and $p$ is the fluid pressure. It is well-known (Langlois \cite{Langlois}) that an incompressible solution of the Stokes equations for the velocity field $(u,v)$ can be written in terms of a stream function $\psi(x,y)$ with 
\begin{equation}
u = {\partial \psi \over \partial y}, \qquad v = - {\partial \psi \over \partial x}.
\end{equation}
The stream function satisfies the biharmonic equation
\begin{equation}
\nabla^4 \psi = 0,
\label{eq:eq0}
\end{equation}
where $\nabla^2$ is the two-dimensional Laplacian. The general solution of (\ref{eq:eq0}) can be written
\begin{equation}
\psi = {\rm Im}[\overline{z} f(z) + g(z)],
\end{equation}
where $f(z)$ and $g(z)$ are analytic functions (which can have isolated singularities) in the fluid region and are often referred to as {\em Goursat functions} (Langlois \cite{Langlois}). These analytic functions are related to physical quantities via 
\begin{equation}
\begin{split}
4 f'(z) &= {p \over \eta} - {\rm i} \omega, \\
-\overline{f(z)} + \overline{z}{f'(z)} +{g'(z)} & = u - {\rm i} v,
\end{split}
\label{eq:eqUV0}
\end{equation}
where $\omega$ is the fluid vorticity.

This connection with analytic functions is important because we propose to harness the powerful mathematical results of analytic function theory to find expressions for the fundamental singularities of Stokes flows which can be evaluated efficiently. Since we will consider singly periodic configurations with
period $l$ in the $x$-direction and channel width $h$ it is natural to non-dimensionalise lengths with respect to the period width $l$, say.

\section{Fundamental singularities of Stokes flow}

In this section we express what is meant by the fundamental singularities of Stokes flow in the language of the isolated singularities of the analytic functions $f(z)$ and $g'(z)$. We focus here on the Stokeslet, stresslet, force quadrupole, along with some irrotational (or ``source'') singularities. Higher-order singularities can be treated similarly.

\subsection{Free-space Stokeslet}

A Stokeslet at a point $z_0$ is associated with a singularity of $f(z)$ having the form
\begin{equation}
f(z) = {\mu \log(z-z_0)},
\label{fstokeslet}
\end{equation}
where $\mu$ is a  generally complex constant which sets the timescale of the motion; we will see later that it is related to the strength of the Stokeslet. The complex velocity is given by
\begin{equation}
u - {\rm i} v = - \overline{f(z)} + \overline{z}{f'(z)} +{g'(z)} =- \overline{\mu} \log(\overline{z-z_0}) + {{\mu} \overline{z} \over{z-z_0}} +{g'(z)}.
\end{equation}
In order that the velocity is {both} single-valued {and} logarithmically singular at $z_0$ 
we must pick
\begin{equation}
{g'(z)} = -\overline{\mu} \log({z-z_0}) - {{\mu} \overline{z_0}\over{z-z_0}}
\label{eq:eq1}
\end{equation}
so that
\begin{equation}
u - {\rm i} v = - \overline{\mu} \log |z-z_0|^2 +{\mu}\left({\overline{z-z_0} \over z-z_0 }\right).
\end{equation}
This expression is single-valued and logarithmically singular at $z_0$.
In summary, if  $f(z)$ and $g'(z)$ locally have the form, near $z_0$, 
\begin{equation}
\begin{split}
f(z) &= {\mu \log(z-z_0)} + {\rm analytic~function},
\\ 
{g'(z)} &= -\overline{\mu} \log({z-z_0}) - {{\mu} \overline{z_0} \over{z-z_0}} + {\rm analytic~function},
\end{split}
\label{eq:eqStokeslet}
\end{equation}
where $\mu \in \mathbb{C}$, then we say there is a Stokeslet at $z_0$.

We can introduce the stress tensor
$\sigma_{ij} = - p \delta_{ij} + 2 \eta e_{ij}$ where $e_{ij}$ is the usual
fluid rate-of-strain tensor.
If we adopt the usual convention of writing a two-dimensional
vector quantity ${\bf a} = (a_x, a_y)$ in complex form as $a_x+ {\rm i} a_y$ then
it can be shown \cite{Langlois} that the complex form of the fluid stress $\sigma_{ij} n_j$ exerted across some point on a contour in the fluid with outward normal $n_i$ is given by
\begin{equation}
2 \eta {\rm i} {\mathrm{d}H \over \mathrm{d}s}, \qquad
H \equiv f(z) + z \overline{f'(z)} + \overline{g'(z)}.
\end{equation}
The force exerted on the fluid by a Stokeslet at  $z_0$ is therefore given by
\begin{equation}
\lim_{\epsilon \to 0}
\oint_{|z-z_0|=\epsilon} 2 \eta {\rm i} {\mathrm{d}H \over \mathrm{d}s} \mathrm{d}s 
= \biggl [  2 \eta {\rm i} H \biggr ]_{|z-z_0|=\epsilon},
\end{equation}
where the square brackets denote the change in the quantity they enclose on a single anticlockwise traversal of the circle ${|z-z_0|=\epsilon}$. This equals
\begin{equation}
\biggl [  2 \eta {\rm i} \left ( \mu \log (z-z_0) - \mu \log (\overline{z-z_0} ) \right ) \biggr ]_{|z-z_0|=\epsilon}
= - 8 \pi \eta \mu,
\label{Strength}
\end{equation}
and gives what is usually known as the strength of the Stokeslet.

\subsection{Free-space stresslet}

Suppose that $f(z)$ has a simple pole at $z_{0}$ so that
\begin{equation}
f(z)=\frac{\mu}{z-z_{0}}, \label{f stresslet}
\end{equation}
where $\mu \in \mathbb{C}$. The associated fluid velocity is encoded in the complex quantity 
\begin{equation}
u - {\rm i} v = - {\overline{\mu} \over \overline{z-z_0}} - {{\mu \overline{z}} \over{(z-z_0})^{2}}+{g'(z)}. 
\end{equation}
In order that this complex velocity field is singular like $1/|z-z_0|$, we must choose
\begin{equation}
g'(z) = {\mu \overline{z_0} \over (z-z_0)^2}, \label{g' stresslet}
\end{equation}
so that
\begin{equation}
u - {\rm i} v = - {\overline{\mu} \over \overline{z-z_0}} - {{\mu} \over{z-z_0}} \left ({\overline{z-z_0} \over z-z_0 }
\right ).
\end{equation}
In summary, if  $f(z)$ and $g'(z)$ locally have the form, near $z_0$, 
\begin{equation}
\begin{split}
f(z) &= {\mu \over z-z_0} +
{\rm analytic~function},\\ 
{g'(z)} &= {\mu \overline{z_0} \over (z-z_0)^2} +{\rm analytic~function},
\end{split}
\label{eq:eqstresslet}
\end{equation}
where $\mu \in \mathbb{C}$, then we say there is a stresslet at $z_0$. The parameter $\mu$ determines the strength of the stresslet.

 \subsection{Free-space force quadrupole}
 
 Suppose that $f(z)$ has a double pole at $z_{0}$ so that
 \begin{equation}
 f(z) = {\mu \over (z-z_0)^2},
\end{equation}
where $\mu \in \mathbb{C}$. Then the associated velocity field can be written as
\begin{equation}
u - {\rm i} v = - {\overline{\mu} \over (\overline{z-z_0})^2} - {2 \mu \overline{z} \over (z-z_0)^3}+ g'(z).
\end{equation}
In order that this complex velocity field is singular like $1/|z-z_0|^{2}$, we must choose
\begin{equation}
g'(z) = {2 \mu \overline{z_0} \over (z-z_0)^3},
\end{equation}
which implies that
\begin{equation}
u - {\rm i} v = -{\overline{\mu} \over (\overline{z-z_0})^2}
- {2 \mu (\overline{z-z_0}) \over (z-z_0)^3}.
\end{equation}
In summary, if  $f(z)$ and $g'(z)$ locally have the form, near $z_0$, 
\begin{equation}
\begin{split}
f(z) &= {\mu \over (z-z_0)^2} + {\rm analytic~function}, \\
g'(z) &= {2 \mu \overline{z_0} \over (z-z_0)^3}+ {\rm analytic~function}, 
\end{split}
\label{eq:eqforcequadrupole}
\end{equation}
where $\mu \in \mathbb{C}$, then we say there is a force quadrupole at $z_0$. The parameter $\mu$ determines the strength of the quadrupole.

\subsection{Source singularities}

The function $g'(z)$ can additionally have its own singularities that are independent of those 
of $f(z)$; we shall refer to these as source singularities. Suppose $f(z)$ is analytic and $g(z)$ has a logarithmic pole at $z_{0}$, i.e.,
\begin{equation}
g(z)=\mu \log(z-z_{0}).
\end{equation}
The Goursat functions have, therefore, the local form near $z_0$ given by
\begin{equation}
\begin{split}
f(z)&={\rm analytic~function}, \\
g'(z)&=\frac{\mu}{z-z_{0}}+ {\rm analytic~function}.
\end{split}
\end{equation} 
If $\mu \in \mathbb{R}$ we say that there is a source/sink at $z_{0}$ (with an associated
mass flux given by $2 \pi \mu$). If $\mu \in {\rm i} \mathbb{R}$ we say that there is a rotlet at $z_{0}$
(which exerts a torque on the fluid of strength dictated by the modulus of $\mu$).

If $g(z)$ has a simple pole at $z_{0}$, i.e.,
\begin{equation}
g(z)=\frac{\mu}{z-z_{0}},
\end{equation}
where $\mu \in \mathbb{C}$, then the Goursat functions
\begin{equation}
\begin{split}
f(z)&={\rm analytic~function}, \\
g'(z)&=-\frac{\mu}{(z-z_{0})^{2}}+ {\rm analytic~function},
\end{split}
\label{sourcedipole}
\end{equation} 
are for a (source) dipole at $z_{0}$. The parameter $\mu$ determines the strength of the dipole.

If $g(z)$ has now a second-order pole at $z_{0}$, i.e.,
\begin{equation}
g(z)=\frac{\mu}{(z-z_{0})^{2}},
\end{equation}
where $\mu \in \mathbb{C}$, then the Goursat functions
\begin{equation}
\begin{split}
f(z)&={\rm analytic~function}, \\
g'(z)&=-\frac{2 \mu}{(z-z_{0})^{3}}+ {\rm analytic~function},
\end{split}
\label{sourcequadrupole}
\end{equation} 
are for a (source) quadrupole at $z_{0}$. The parameter $\mu$ determines the strength of the quadrupole.

\section{Conformal mapping}

\begin{figure}
\begin{center}
\includegraphics[scale=1]{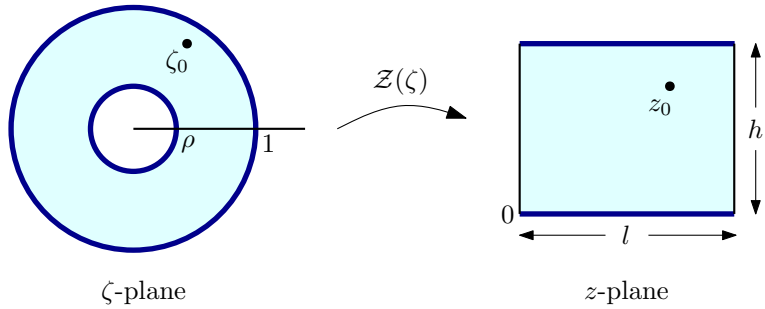}
\caption{Conformal mapping from the annulus $\rho < |\zeta| < 1$ in $\zeta$-plane to the period rectangle in the physical $z$-plane occupying the region $0 \le x \le 2 \pi=l$, $0 \le y \le -\log \rho =h$.  \label{Fig1}}
\end{center}
\end{figure}

To study singly periodic arrays of singularities, we will exploit the conformal mapping
\begin{equation}
z= {\cal Z}(\zeta) \equiv - {\rm i} \log \zeta.
\label{eq:eqmap}
\end{equation}
This map transplants 
the annulus $\rho < |\zeta| < 1$ to a period rectangle in the $z=x+{\rm i} y$
plane occupying the region
\begin{equation}
0 \le x \le 2 \pi=l, \qquad 0 \le y \le -\log \rho=h,
\end{equation}
as illustrated in Fig. \ref{Fig1}. Changing the $x$-period simply requires multiplication of (\ref{eq:eqmap}) by the appropriate real factor.  Under this conformal map, the unit circle $|\zeta|=1$ corresponds to $0\leq x \leq l$, $y=0$, the inner circle $|\zeta|=\rho$ corresponds to $0\leq x \leq l$, $y=h$ and the interval $[\rho,1]$ (in $\zeta$-plane) to the vertical sides of the period window $x=0,l$, $0\leq y\leq h$.

The key idea of our approach is to show that analytical expressions, with fast convergence properties, for the fundamental singularities of periodic Stokes flow can be derived in terms of the variable $\zeta$. It is then a simple matter, if required, to re-express the final results as functions of $z$ using the relation
\begin{equation}
\zeta = {\rm e}^{{\rm i} z}
\label{zetaz}
\end{equation}
which follows from (\ref{eq:eqmap}).

The representative point singularity in the fundamental period rectangle (Fig. \ref{Fig1}) is located at $z_0$ and has a  preimage at $\zeta_0$ such that
\begin{equation}
z_{0}={\cal Z}(\zeta_{0})=-{\rm i} \log \zeta_{0}.
\label{zeta0eqn}
\end{equation}

\section{Periodic Stokes singularity arrays in a channel \label{Sec5}}

\begin{figure}
\begin{center}
\includegraphics[scale=1]{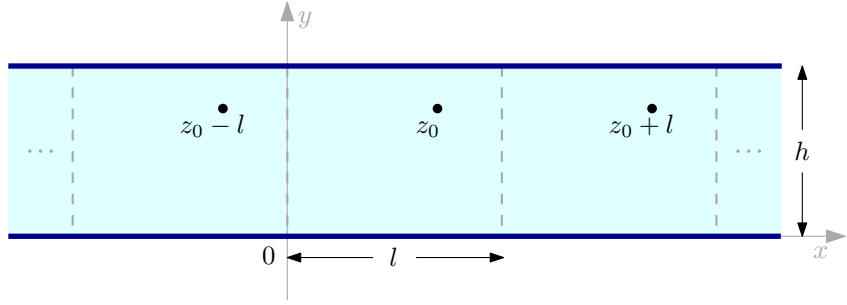}
\caption{A periodic array of point singularities at points  $z=z_{0}+n l$, $n \in \mathbb{Z}$, with $0<\text{Re}[z_{0}]<l$, $0<\text{Im}[z_{0}]<h$, in a two-dimensional channel $-\infty<x<\infty$, $0 \leq y \leq h$.}
\label{fig: configurationchannel}
\end{center}
\end{figure}

Consider a two-dimensional channel $-\infty<x<\infty$, $0 \leq y \leq h$ and a periodic array of point singularities at $z=z_{0}+n l$, $n \in \mathbb{Z}$, with $0<\text{Re}[z_{0}]<l$, $0<\text{Im}[z_{0}]<h$. Figure \ref{fig: configurationchannel} shows a schematic of the configuration. This problem can be solved using standard Fourier transform techniques (Pozrikidis \cite{Poz}), but we will give an alternative derivation, and form, of the solution by extending the ideas given in the previous sections.

\subsection{Periodic array of Stokeslets in a channel}

Consider a periodic array of Stokeslets each of strength given by \eqref{Strength} for some $\mu \in \mathbb{C}$. The representative Stokeslet in the fundamental period rectangle shown in Fig. \ref{Fig1} is located at $z_{0}$ and has a preimage at $\zeta_{0}$ related to $z_0$ 
via (\ref{zeta0eqn}).
%with
%\begin{equation}
%z_{0}= {\cal Z}(\zeta_{0})=-{\rm i} \log \zeta_{0}.
%\end{equation}
An important step is to introduce the functions
\begin{equation}
F(\zeta) \equiv f({\cal Z}(\zeta)), \qquad G(\zeta) \equiv g'({\cal Z}(\zeta)).
\label{FGzeta}
\end{equation}
Now let
\begin{equation}\label{expansionchannel}
\begin{split}
F(\zeta) &= a \log \zeta + \mu \log (\zeta-\zeta_0) + \hat F(\zeta), \\
G(\zeta) &= - \overline{a} \log \zeta - \overline{\mu} \log (\zeta-\zeta_0) + {\lambda \over \zeta-\zeta_0}
+ {\rm i} \log \zeta \left [ {F'(\zeta) \over {\cal Z}'(\zeta)} \right ] + \hat G(\zeta),
\end{split}
\end{equation}
where $\hat F(\zeta)$ and $\hat G(\zeta)$ are analytic and single-valued in the annulus $\rho < |\zeta| < 1$ and the constants $\lambda$ and $a$ are to be found.
With this ansatz, we have
\begin{equation}
{F'(\zeta) \over {\cal Z}'(\zeta)} = {\rm i} \zeta \left [ {a \over \zeta} + {\mu \over \zeta-\zeta_0} + \hat F'(\zeta) \right ] =
{\rm i} a + {{\rm i} \mu \zeta \over \zeta-\zeta_0} + {\rm i} \zeta \hat F'(\zeta)
\end{equation}
and, hence, the complex velocity field is
\begin{equation}
\begin{split}
u - {\rm i} v = - \overline{a} \log |\zeta|^2 
&- \overline{\mu} \log |\zeta-\zeta_0|^2
+ {\rm i} \log |\zeta|^2 \left [ {F'(\zeta) \over {\cal Z}'(\zeta)} \right ]
\\ &- \overline{\hat F(\zeta)} + \hat G(\zeta) + {\lambda \over \zeta-\zeta_0}.
\end{split}
\end{equation}
By the choice of the ansatz, $u - {\rm i} v$ is invariant as $\zeta \mapsto \zeta {\mathrm e}^{2 \pi {\rm i}}$ and this corresponds to the periodicity of the velocity field along the channel. A local analysis of the singularity of $u-{\rm i}v$ necessitates that we choose
\begin{equation}
\lambda = \mu \zeta_0 \log |\zeta_0|^2,
\end{equation}
in order to ensure that the singularity at $\zeta_0$ has the required form.

The no-slip condition on $|\zeta|=1$ becomes
\begin{equation}
- \overline{\hat F(\zeta)} + \hat G(\zeta) = \overline{\mu}  \log |\zeta-\zeta_0|^2 - {\mu \zeta_0 \log |\zeta_0|^2 \over \zeta-\zeta_0},
\label{eq:eqBC1channel}
\end{equation}
while the no-slip condition on $|\zeta|=\rho$ is
\begin{equation}
\begin{split}
- \overline{\hat F(\zeta)} + \hat G(\zeta) 
&- \log \rho^2 [ \zeta \hat F'(\zeta) ] - (a+ \overline{a}) \log \rho^2 
\\&= \overline{\mu}  \log |\zeta-\zeta_0|^2  - {\mu \zeta_0 \log |\zeta_0|^2 \over \zeta-\zeta_0} 
+ { \mu \zeta  \log \rho^2\over \zeta-\zeta_0}.
\end{split}
\label{eq:eqBC2channel}
\end{equation}
For $|\zeta|=1$, conditions (\ref{eq:eqBC1channel}) and (\ref{eq:eqBC2channel}) can be written in the form
\begin{equation}
\begin{split}
- \overline{\hat F}(\zeta^{-1}) + \hat G(\zeta) &= \sum_{n=-\infty}^\infty d_n \zeta^n, \\
- \overline{\hat F}(\rho \zeta^{-1}) + \hat G(\rho \zeta) 
- \log \rho^2 [\rho \zeta \hat F'(\rho \zeta) ] - (a+ \overline{a}) \log \rho^2 
&=  \sum_{n=-\infty}^\infty e_n \zeta^n,
\end{split}
\label{eq:eqBC4channel}
\end{equation}
where we have introduced the following Laurent expansions of known functions:
\begin{equation}
\begin{split}
\sum_{n=-\infty}^\infty d_n \zeta^n &=
\overline{\mu}  \log |\zeta-\zeta_0|^2 - {\mu \zeta_0 \log |\zeta_0|^2 \over \zeta-\zeta_0}, \\
\sum_{n=-\infty}^\infty e_n \zeta^n &=
\overline{\mu}  \log |\rho \zeta-\zeta_0|^2  - {\mu \zeta_0 \log |\zeta_0|^2 \over \rho \zeta-\zeta_0} 
+ { \mu\rho \zeta  \log \rho^2\over \rho \zeta-\zeta_0}.
\end{split}
\end{equation}
The coefficients $\{ d_{n}, e_{n} | n \in \mathbb{Z} \}$ can be determined explicitly by computing the relevant Laurent expansions. Alternatively, they can be computed numerically using Fast Fourier transforms. Using local expansions, we find
\begin{equation}
\begin{split}
d_{n}&=-\overline{\mu} (\overline{\zeta_{0}})^{n}/n, \quad n \ge 1, \\
d_{-n}&=-(\zeta_{0})^{n}(\overline{\mu}/n+\mu \log|\zeta_{0}|^{2}), \quad n \ge 1, \\
d_{0}&=0
\end{split}\label{dnStokeslet}
\end{equation}
and
\begin{equation}
\begin{split} 
e_{n}&=\rho^{n} (\zeta_{0})^{-n}[\mu (\log|\zeta_{0}|^{2}-\log \rho^{2})-\overline{\mu}/n], \quad n \ge 1, \\
e_{-n}&=-\overline{\mu} \rho^{n} (\overline{\zeta_{0}})^{-n}/n, \quad n \ge 1, \\
e_{0}&=2 \text{Re}[\mu] \log|\zeta_{0}|^{2}.
\end{split}\label{enStokeslet}
\end{equation}

The next step is to consider the Laurent series expansions
\begin{equation}\label{Laurentexpchannel}
\begin{split}
\hat F(\zeta) &= \sum_{n=1}^\infty F_n \zeta^n + \sum_{n=1}^\infty H_n \left ({\rho \over \zeta} \right )^n, \\
\hat G(\zeta) &= \sum_{n=1}^\infty G_n \zeta^n + G_0 + \sum_{n=1}^\infty K_n \left ({\rho \over \zeta} \right )^n.
\end{split}
\end{equation}
Without loss of generality, the constant term in the expansion of $\hat F(\zeta)$ can be set equal to zero
owing to an additive degree of freedom in the specification of $f(z)$.
Remarkably, the Laurent series (\ref{Laurentexpchannel})  can be substituted into (\ref{eq:eqBC4channel}) and their unknown coefficients found explicitly. The constant terms in (\ref{eq:eqBC4channel}) give
\begin{equation}\label{G0}
G_0 = d_0, \qquad G_0 - \log \rho^2 (a+ \overline{a}) = e_0
\end{equation}
which implies that
\begin{equation}
{\rm Re}[a] = {d_0 - e_0 \over 4 \log \rho}.
\end{equation}
The imaginary part of $a$ can be set to zero because making the transformation $f(z) \mapsto f(z) + b z$ for some $b \in \mathbb{R}$ clearly does not affect the velocity field and this corresponds to adding the term $-{\rm i} b \log \zeta$ to $F(\zeta)$. It follows that we can take
\begin{equation}\label{parameter a}
a= {d_0 - e_0 \over 4 \log \rho}.
\end{equation}
Equating coefficients of the other powers of $\zeta$ produces the equations
\begin{equation}
\begin{split}
-\rho^n \overline{H_n} + G_n &= d_{n}, \\
-\overline{F_n} + \rho^n K_n &= d_{-n}, \\
-\overline{H_n} + \rho^n G_n - n \rho^n \log \rho^2 F_n &= e_{n}, \\
-\rho^n \overline{F_n} + K_n + n \log \rho^2 H_n &= e_{-n}, \quad n \ge 1.
\end{split}
\label{eqnFcoeff}
\end{equation}
This system can be manipulated to find
\begin{equation}\label{coef1channel}
F_n = {n \rho^n \log \rho^2 e_n + \rho^n (1-\rho^{2n}) \overline{e_{-n}} - 
n \rho^{2n} \log \rho^2 d_n - (1-\rho^{2n}) \overline{d_{-n}} \over
(1-\rho^{2n})^2 - n^2 \rho^{2n} (\log \rho^2)^2}, \quad n \ge 1.
\end{equation}
The remaining coefficients needed to evaluate $\hat F(\zeta)$ and $\hat G(\zeta)$ then follow
 by back substitution:
\begin{equation}\label{coef2channel}
G_{n} = {-n \rho^{2n} \log \rho^2 F_n - \rho^n e_n + d_n \over (1-\rho^{2n})}, \quad
H_{n} = \frac{\overline{G_n} - \overline{d_n}}{\rho^n}, \quad 
K_{n} = \frac{\overline{F_n} + d_{-n}}{\rho^n}, ~ n \ge 1.
\end{equation}
It should be clear that substitution of the explicit form (\ref{eqnFcoeff}) into 
(\ref{coef2channel}) will lead to explicit expressions for all these coefficients
although we have avoided displaying those formulas here.

\begin{figure}
\begin{center}
\includegraphics[scale=0.25]{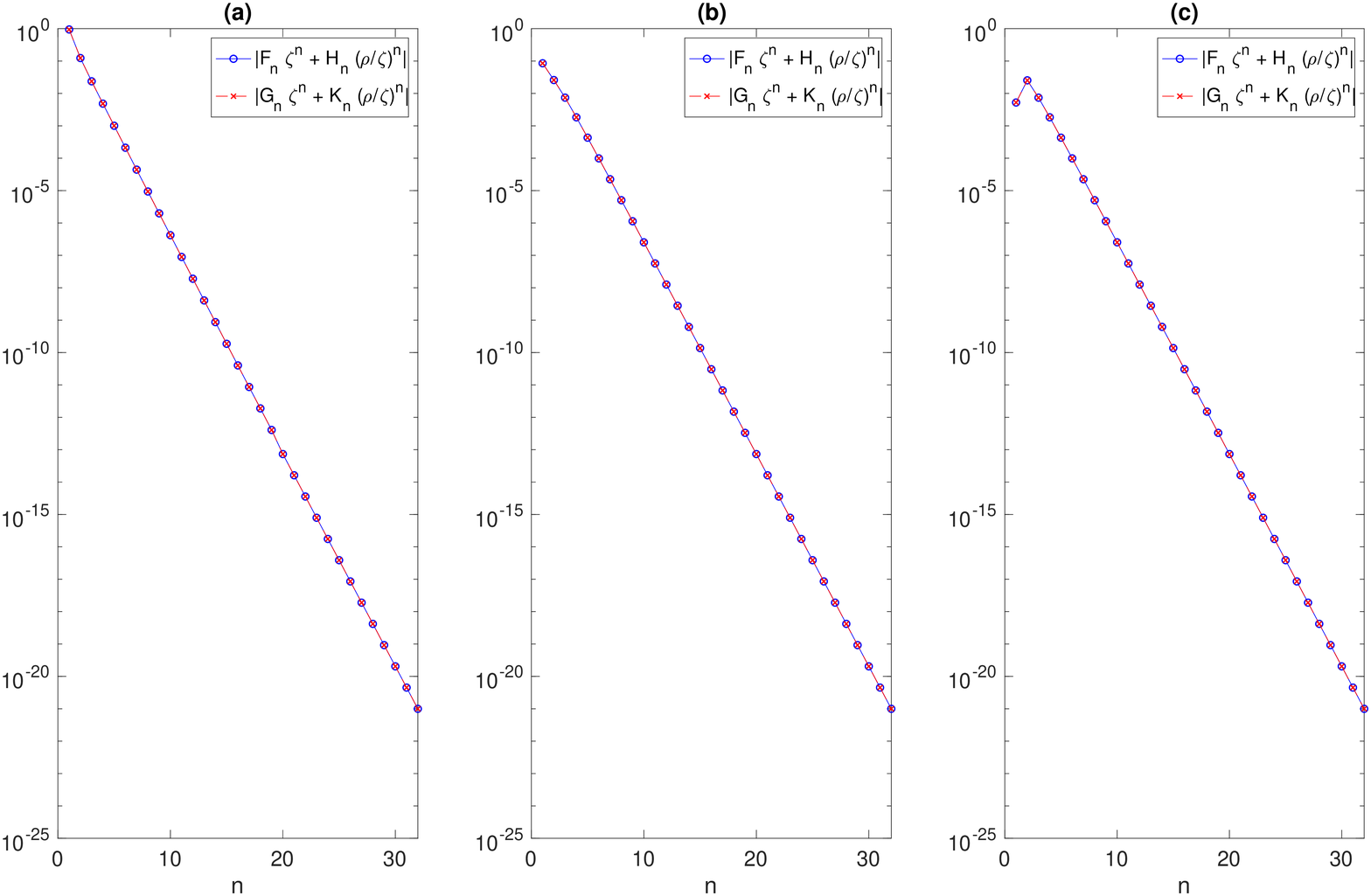}
\caption{Decay of the Laurent expansion terms $F_{n} \zeta^n + H_n (\rho/\zeta)^n$ and $G_{n} \zeta^n + K_n (\rho/\zeta)^n$, $n \geq 1$, given by expressions \eqref{Laurentexpchannel}, for parameters $l=2\pi$, $\zeta=\mathrm{e}^{-1}$, $\zeta_0=0.6$ and different channel heights: (a) $h=\pi/2$, (b) $h=\pi$ and (c) $h=2\pi$.}
\label{fig: graph1}
\end{center}
\end{figure}

\begin{figure}
\begin{center}
\includegraphics[scale=0.25]{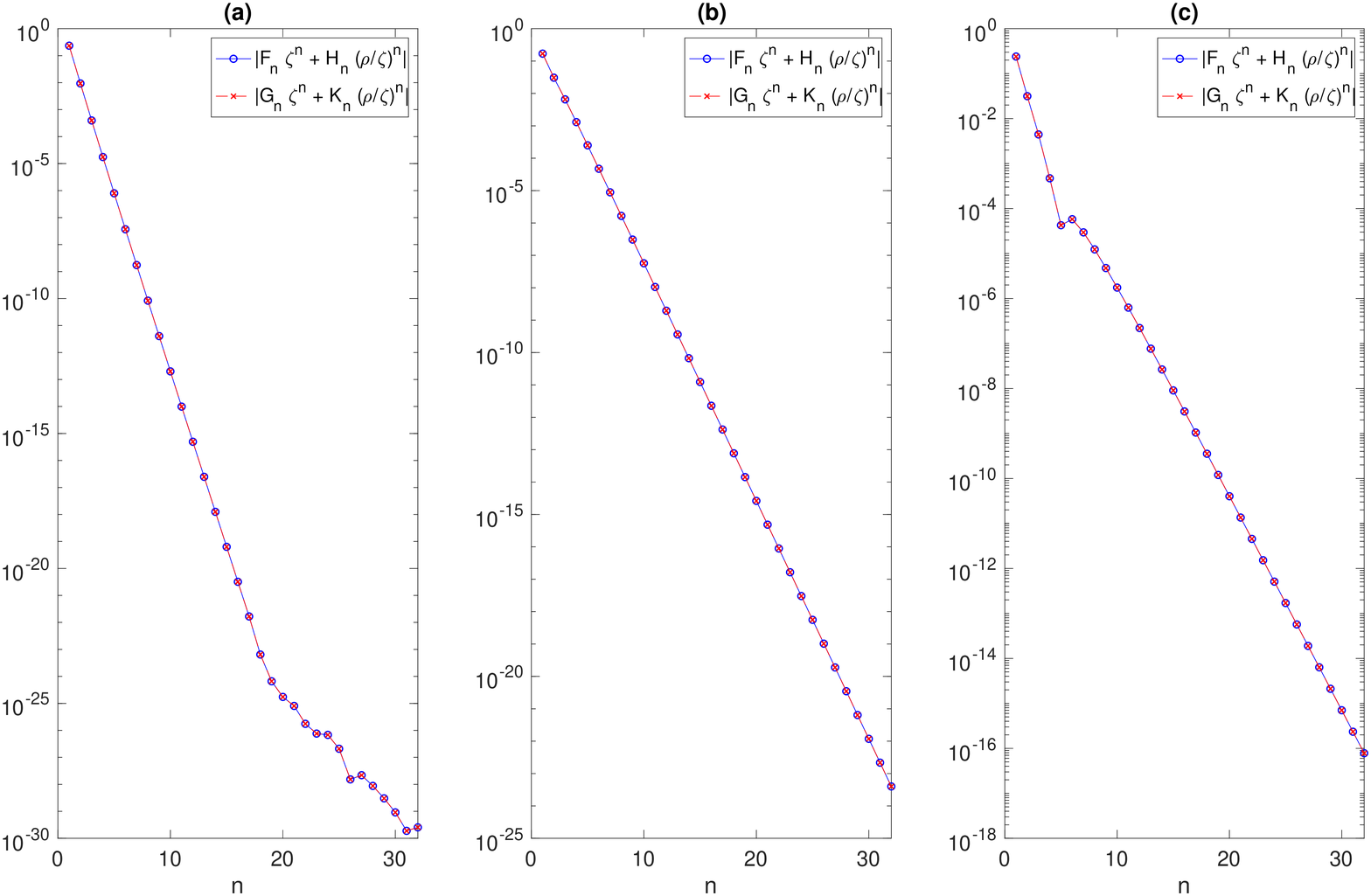}
\caption{Decay of the Laurent expansion terms $F_{n} \zeta^n + H_n (\rho/\zeta)^n$ and $G_{n} \zeta^n + K_n (\rho/\zeta)^n$, $n \geq 1$, given by expressions \eqref{Laurentexpchannel}, for parameters $l=2\pi$, $h=\pi$, $\zeta=\mathrm{e}^{-1}$ and different locations for the principal point Stokeslet: (a) $\zeta_0=0.1$, (b) $\zeta_0=0.5$ and (c) $\zeta_0=0.9$.}
\label{fig: graph2}
\end{center}
\end{figure}

In summary, the Goursat functions for a periodic array of Stokeslets in a channel with one at $z_{0}$ in the principal period window and having strength given by \eqref{Strength} are
\begin{equation}
\begin{split}
F(\zeta) &= a \log \zeta + \mu \log (\zeta-\zeta_0) + \sum_{n=1}^\infty F_n \zeta^n + \sum_{n=1}^\infty H_n \left ({\rho \over \zeta} \right )^n, \\
G(\zeta) &= - \overline{a} \log \zeta - \overline{\mu} \log (\zeta-\zeta_0) + {\lambda \over \zeta-\zeta_0}
+ {\rm i} \log \zeta \left [ {F'(\zeta) \over {\cal Z}'(\zeta)} \right ] \\
&~~~+ \sum_{n=1}^\infty G_n \zeta^n + G_0 + \sum_{n=1}^\infty K_n \left ({\rho \over \zeta} \right )^n, \\
\text{with} &\qquad a= {d_0 - e_0 \over 4 \log \rho}, \qquad \lambda = \mu \zeta_0 \log |\zeta_0|^2,
\end{split}\label{FG SP Stokeslet}
\end{equation}
where coefficients $\{F_{n}, H_{n}, G_{n}, K_{n}|n \in \mathbb{N}\}$ and $G_{0}$ are given by \eqref{G0}, \eqref{coef1channel} and \eqref{coef2channel} in terms of \eqref{dnStokeslet}--\eqref{enStokeslet}. Figures \ref{fig: graph1} and \ref{fig: graph2} show the fast decay of the Laurent expansion terms \eqref{Laurentexpchannel} and indicate that only a few terms are required in general to provide accurate values of the flow field variables.

Davis \cite{Davis1993} presented solutions to various problems involving distributions of Stokeslets in an unbounded domain, as well as in a two-dimensional channel, with a view to understanding the blocking properties of periodic arrays of wall-attached barriers. We have verified numerically that the associated flow field found above agrees with the solutions found by Davis \cite{Davis1993}.

As a separate check,
yet another
new method for solving the problem of a periodic array of point singularities in a channel geometry is summarized in Appendix \ref{appendixA}. That method is based on a novel transform approach to biharmonic boundary value problems recently described by the authors \cite{LucaCrowdy};
it provides a quasi-analytical solution to the problem
in that it reduces determination of the Goursat functions to the solution of a
small linear system whose coefficients are given by explicit integrals.
  The solutions presented above have been verified numerically against the latter solutions.

%\noindent
%{\bf Davis' solution}
%
%Davis finds the solution for a point force of strength $4 \pi \eta {\bf x}$ where 
%$\eta$ is the fluid viscosity. His velocity field can be written in complex form as
%\begin{equation}
%u - {\rm i} v = - {1 \over 2} \log \overline{z} + {\overline{z} \over z} - {1 \over 2} \log z
%\end{equation}
%so we identify
%\begin{equation}
%f(z) = {1 \over 2} \log z, \qquad g'(z) = - {1 \over 2} \log z
%\end{equation}
%Hence, his solution corresponds to the choice $\mu = 1/2$.
%This is consistent with the fact that the force exerted on the fluid is
%\begin{equation}
%F_x + {\rm i} F_y = 2 \eta {\rm i} \left [ \mu \log \left ({z-z_0 \over \overline{z-z_0}} \right ) \right ]_{|z-z_0|=\epsilon}
%= 2 \eta {\rm i} (-4 \pi {\rm i})  \mu = 8  \pi \eta \mu = 4  \pi \eta.
%\end{equation}

%Figure 2 shows the streamline pattern for a doubly periodic array of Stokeslets of strength (13)with ? = 1 located
%at the centre of each period window. To plot these streamlines, the rapidly convergent representations (46)�(51)
%which appear in (68) have been truncated to only N = 16. Similar plots computed using different methods are
%presented in Pozrikidis [15].

\begin{figure}
  \centering
  \begin{tabular}{@{}c@{}}
    \includegraphics[scale=0.3]{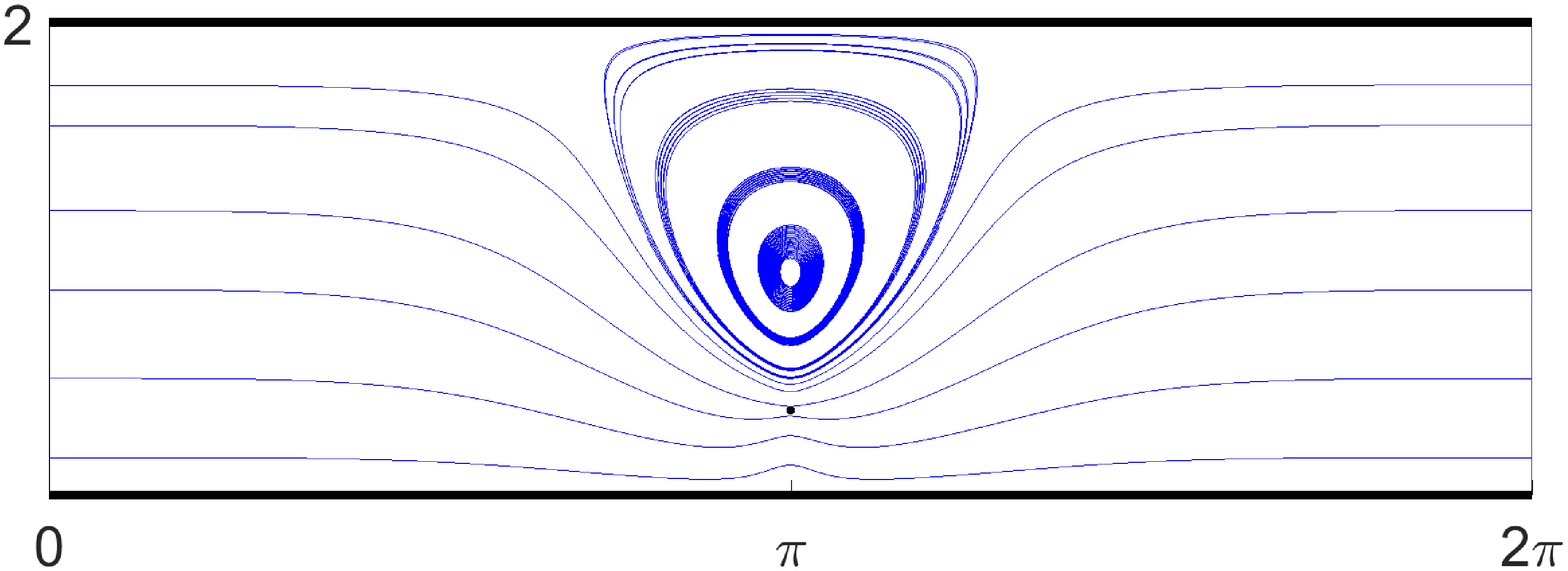} \\[\abovecaptionskip]
  \end{tabular}

\vspace{-0.6in}
  \vspace{\floatsep}

  \begin{tabular}{@{}c@{}}
    \includegraphics[scale=0.3]{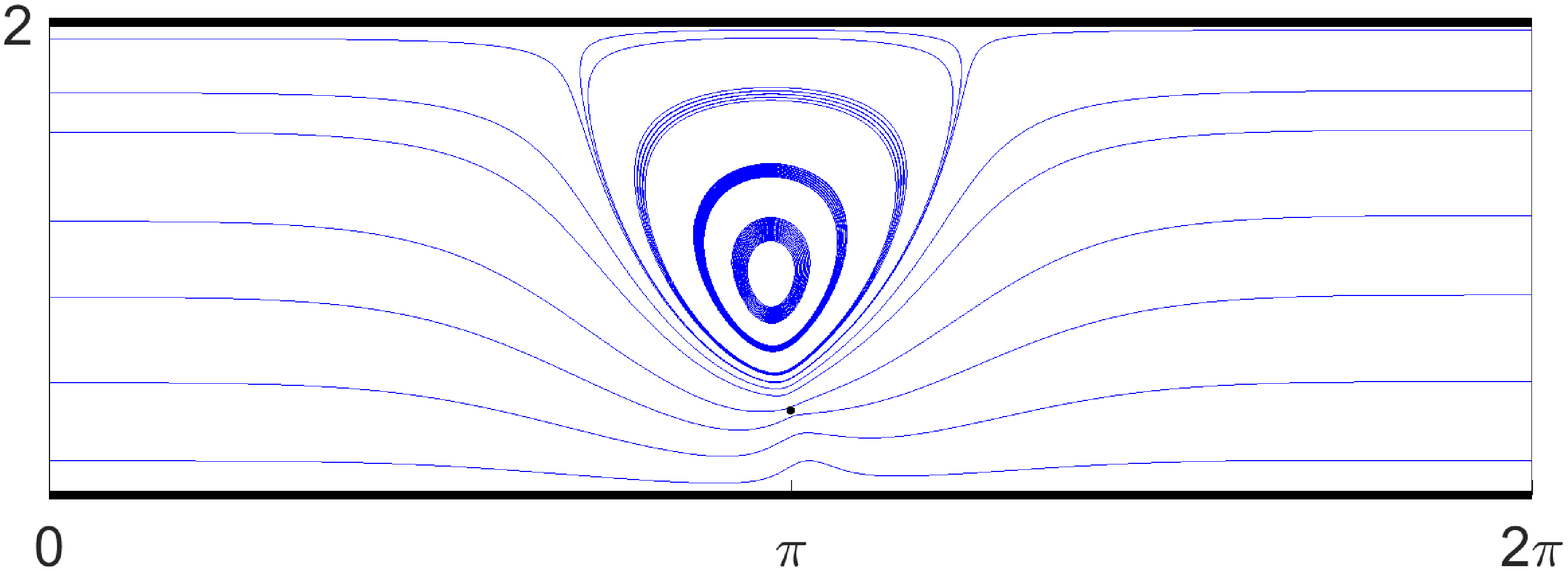} \\[\abovecaptionskip]
  \end{tabular}

\vspace{-0.6in}
  \vspace{\floatsep}

  \begin{tabular}{@{}c@{}}
    \includegraphics[scale=0.3]{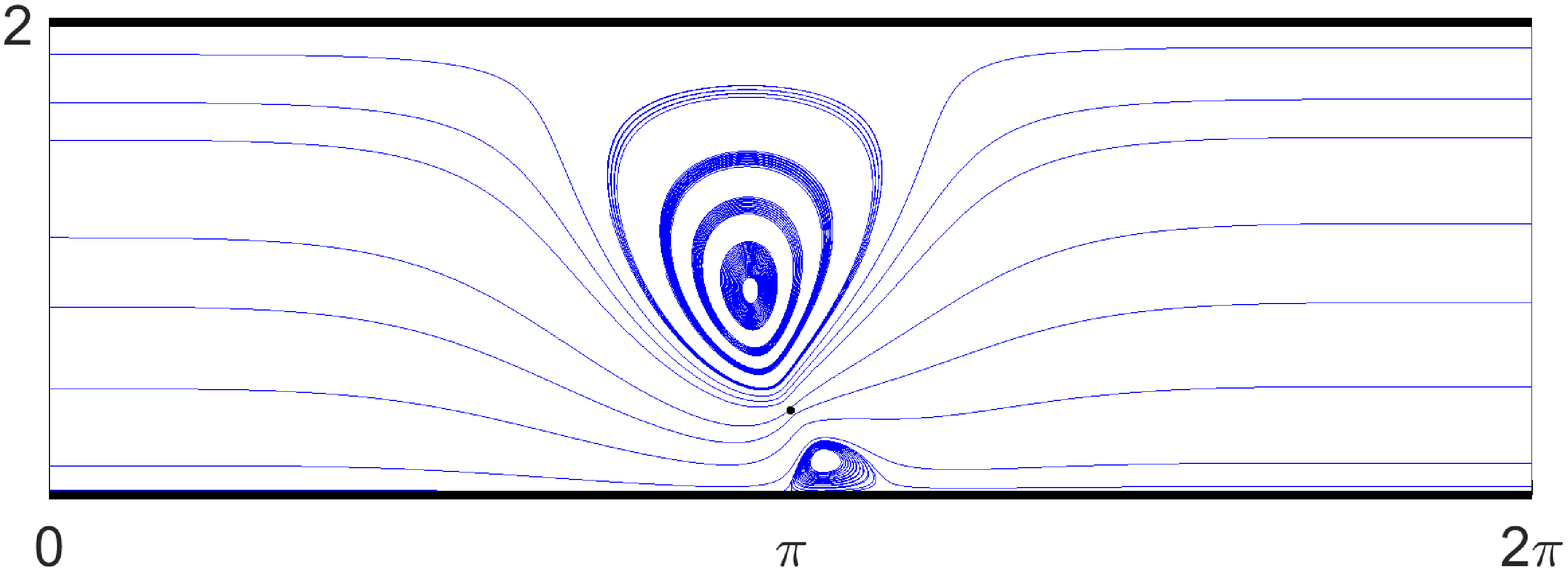} \\[\abovecaptionskip]
  \end{tabular}

\vspace{-0.6in}
 \vspace{\floatsep}

  \begin{tabular}{@{}c@{}}
    \includegraphics[scale=0.3]{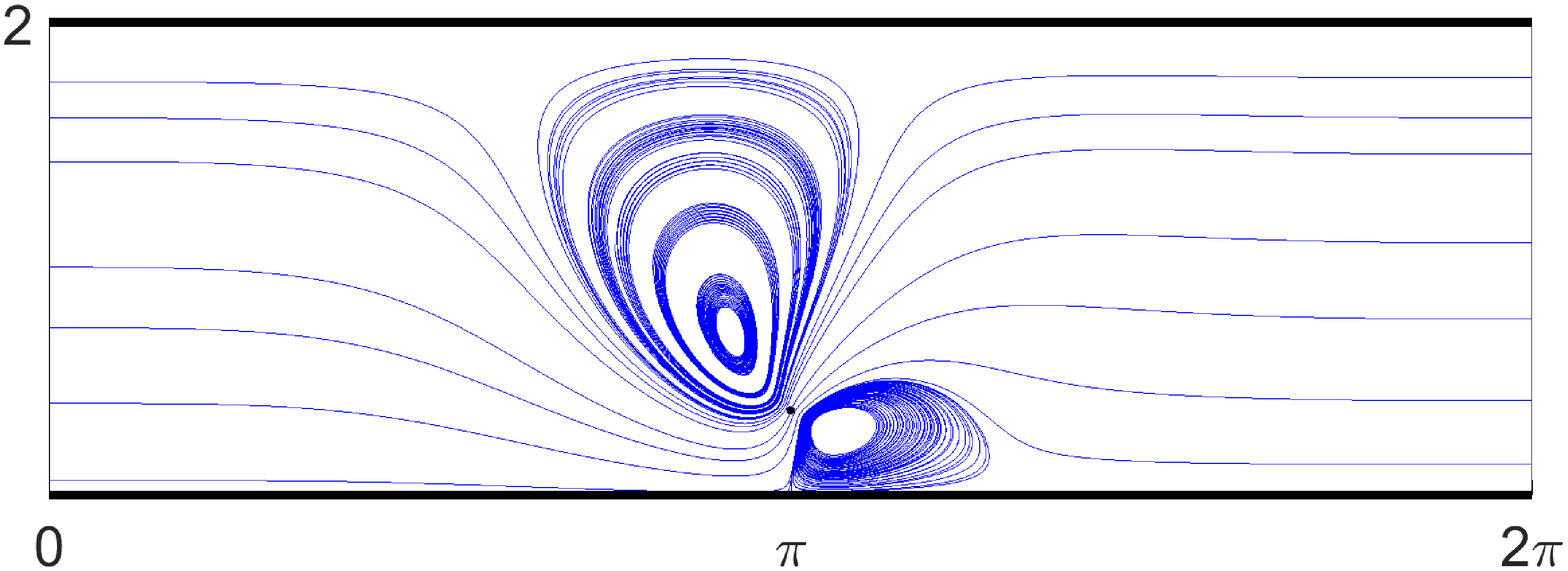} \\[\abovecaptionskip]
  \end{tabular}

\vspace{-0.6in}
 \vspace{\floatsep}

  \begin{tabular}{@{}c@{}}
    \includegraphics[scale=0.3]{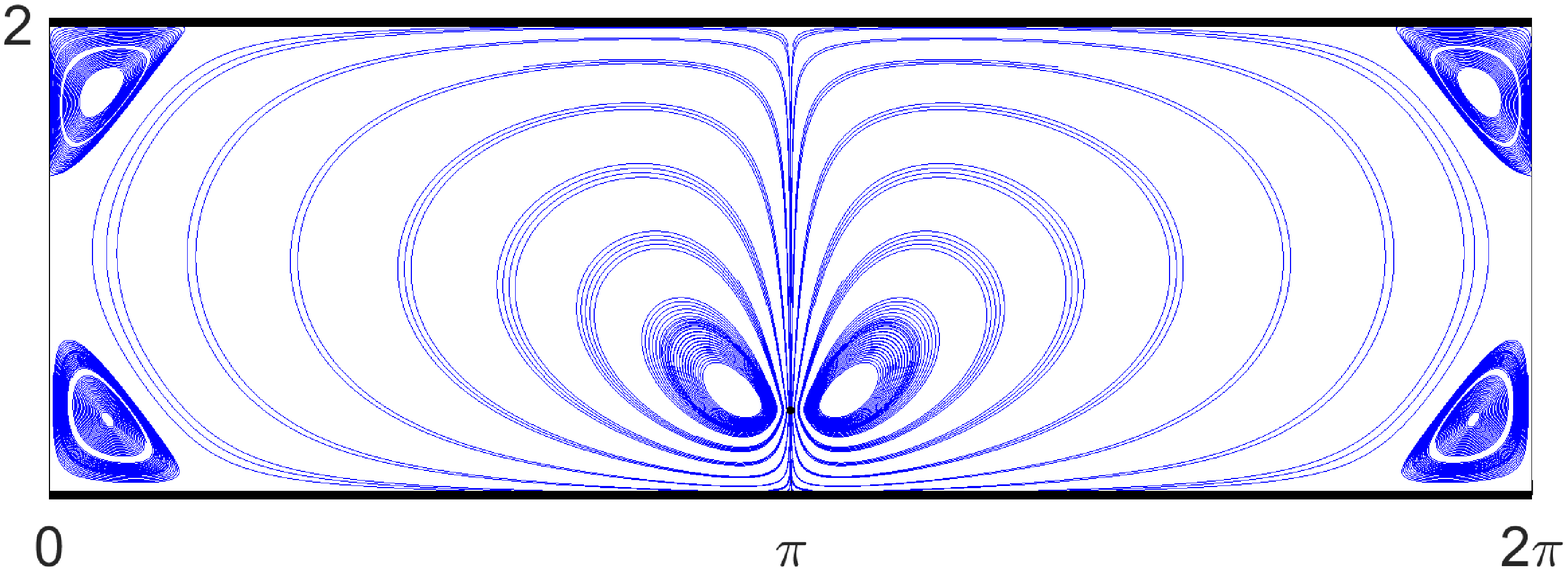} \\[\abovecaptionskip]
  \end{tabular}

\vspace{-0.3in}

\caption{Streamline patterns for a periodic Stokeslet in a channel in a period window for different strengths $\mu$. The principal period window is $0 \leq x \leq l=2 \pi$, $0 \leq y \leq h=2$. The principal point Stokeslet is located at $z_{0}=\pi-{\rm i} \log (0.7)$ and its strength is $\mu=\mathrm{e}^{{\rm i}\phi}$, (a) $\phi=0$, (b) $\phi=\pi/8$, (c) $\phi=\pi/4$, (d) $\phi=3\pi/8$ and (e) $\phi=\pi/2$.}

\label{fig. Stokeslets strengths}
\end{figure}

\begin{figure}
  \centering
  \begin{tabular}{@{}c@{}}
    \includegraphics[scale=0.3]{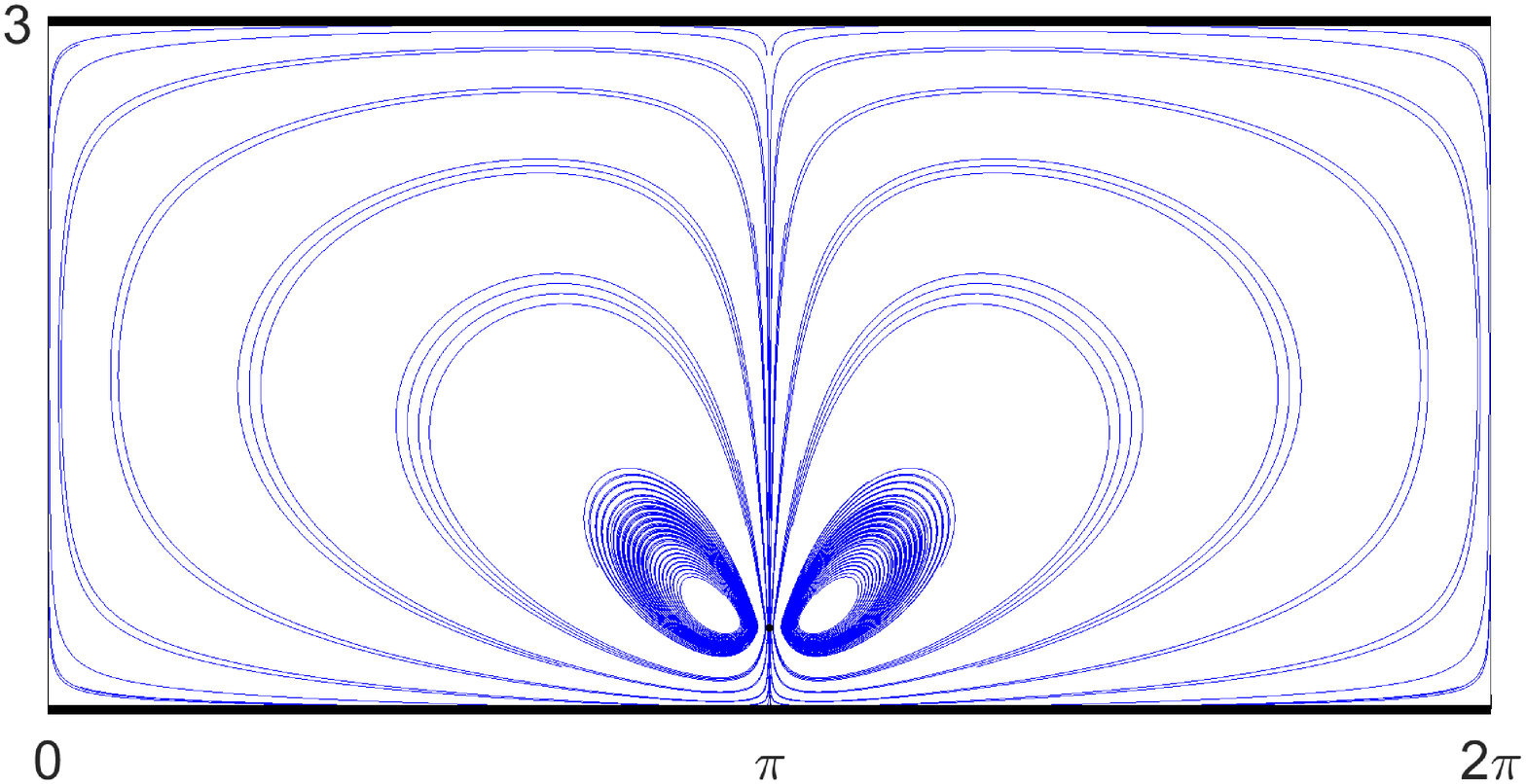} \\[\abovecaptionskip]
  \end{tabular}

\vspace{-0.3in}
  \vspace{\floatsep}

  \begin{tabular}{@{}c@{}}
    \includegraphics[scale=0.3]{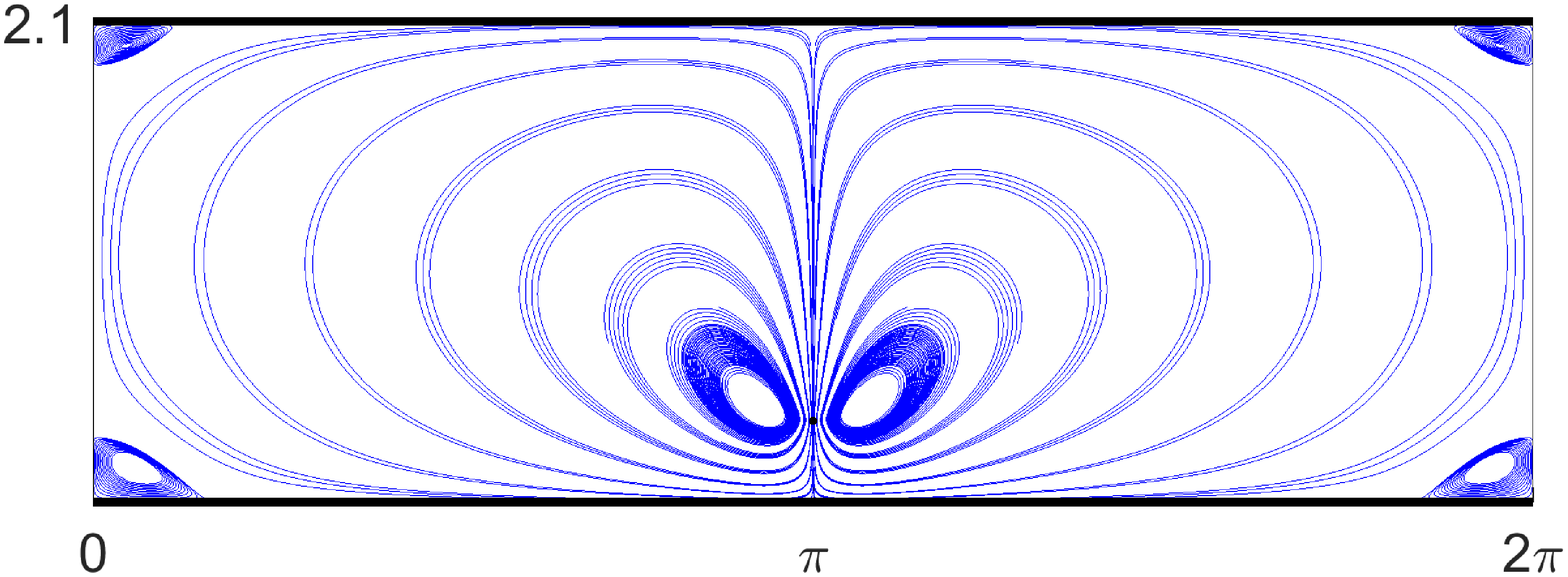} \\[\abovecaptionskip]
  \end{tabular}

\vspace{-0.55in}
  \vspace{\floatsep}

  \begin{tabular}{@{}c@{}}
    \includegraphics[scale=0.3]{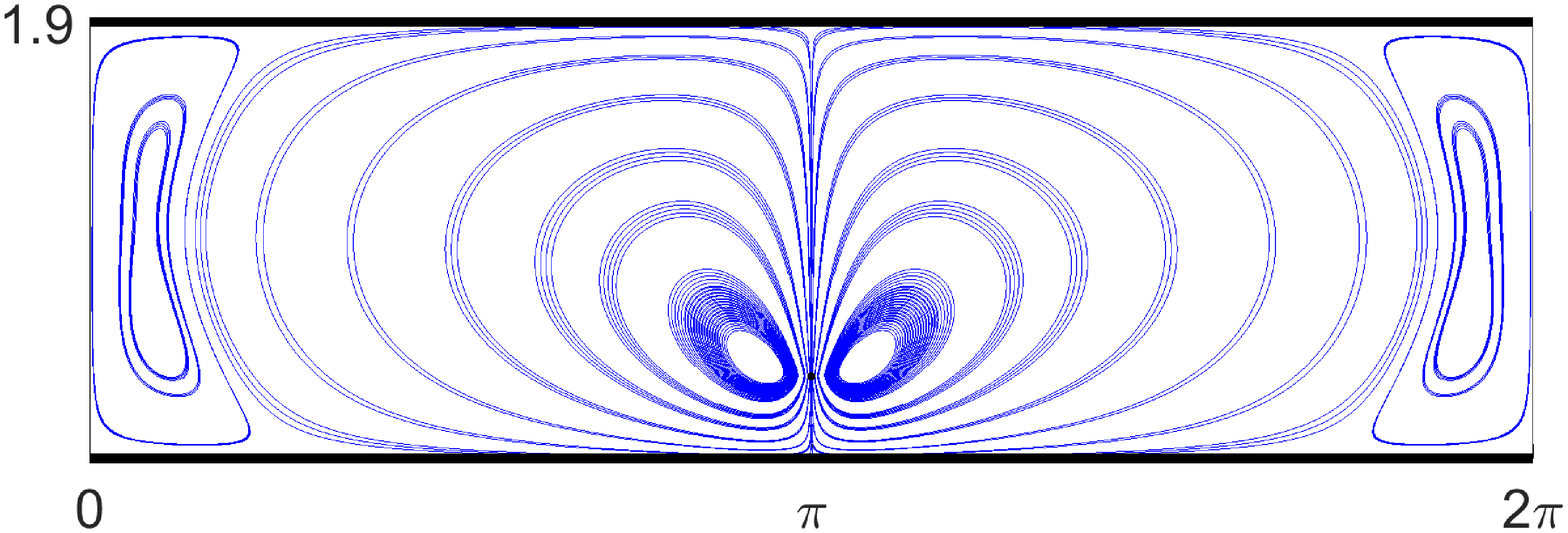} \\[\abovecaptionskip]
  \end{tabular}

\vspace{-0.7in}
  \vspace{\floatsep}

  \begin{tabular}{@{}c@{}}
    \includegraphics[scale=0.3]{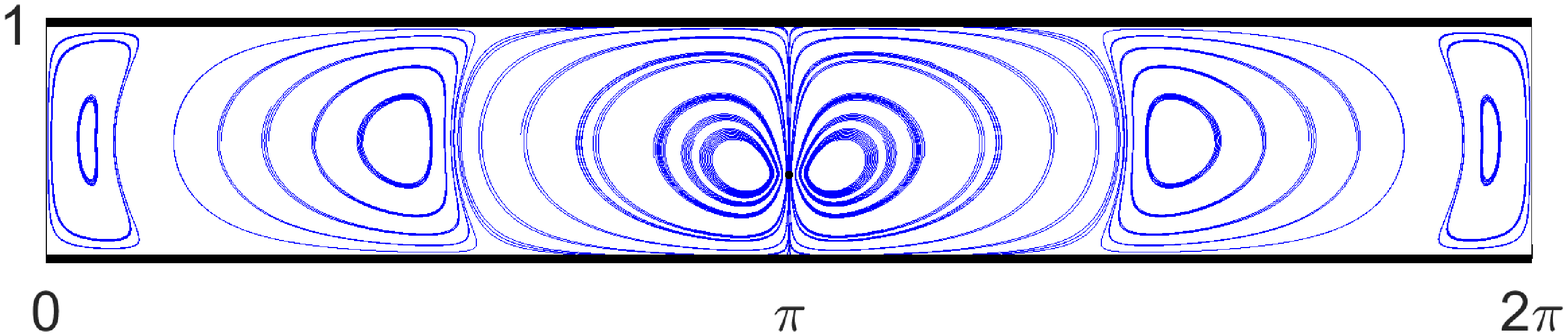} \\[\abovecaptionskip]
  \end{tabular}

\vspace{-0.9in}
  \vspace{\floatsep}

  \begin{tabular}{@{}c@{}}
    \includegraphics[scale=0.3]{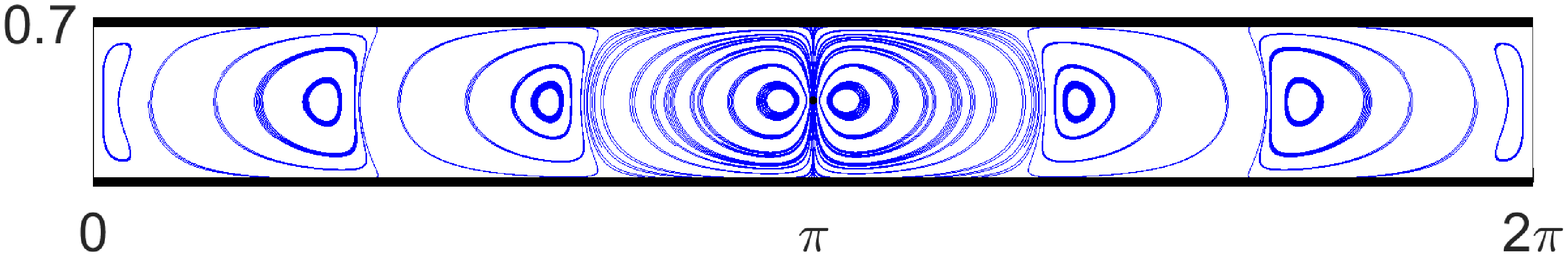} \\[\abovecaptionskip]
  \end{tabular}

\vspace{-0.5in}

\caption{Streamline patterns for a periodic Stokeslet in a channel in a period window for different channel heights $h=0.7, 1, 1.9, 2.1, 3$. The principal period window is $0 \leq x \leq l=2 \pi$, $0 \leq y \leq h$. The principal point Stokeslet is located at $z_{0}=\pi-{\rm i} \log (0.7)$ and its strength is $\mu={\rm i}$.}

\label{fig. Stokeslets heights}
\end{figure}

Figure \ref{fig. Stokeslets strengths} shows the streamline patterns for a singly periodic array of Stokeslets of different strengths $\mu=\mathrm{e}^{{\rm i}\phi}$, $\phi=0, \pi/8, \pi/4, 3\pi/8, \pi/2$ and parameters $l=2\pi$, $h=2$, $z_{0}=\pi-{\rm i} \log (0.7)$ which is the location of the principal point Stokeslet. The graphs show the streamline topology transition between $\phi \rightarrow \pi/2$ and $\phi=\pi/2$ (other singularities behave similarly). Figure \ref{fig. Stokeslets heights} shows the streamline patterns for a singly periodic array of Stokeslets of strength $\mu={\rm i}$ for different channel heights $h=0.7, 1, 1.9, 2.1, 3$ and parameters $l=2\pi$, $z_{0}=\pi-{\rm i} \log (0.7)$.
%(last case of figure \ref{fig. Stokeslets strengths})
%For only this strength choice, we obtain these streamlines patterns.

%-----------------------------------------%

%-----------------------------------------%

%-----------------------------------------%

\subsection{Periodic array of stresslets in a channel}

Following similar steps, we find that the Goursat functions for a periodic array of torque-free stresslets in a channel with local behaviour \eqref{eq:eqstresslet} at $z_{0}$ in the principal period window are
\begin{equation}
\begin{split}
F(\zeta) &= a \log \zeta + \frac{{\rm i} \mu \zeta_{0}}{\zeta-\zeta_0} + \sum_{n=1}^\infty F_n \zeta^n + \sum_{n=1}^\infty H_n \left ({\rho \over \zeta} \right )^n, \\
G(\zeta) &= - \overline{a} \log \zeta + \frac{\chi}{\zeta-\zeta_0} + \frac{\nu}{(\zeta-\zeta_0)^{2}}
+ {\rm i} \log \zeta \left [ {F'(\zeta) \over {\cal Z}'(\zeta)} \right ] \\
&~~~+ \sum_{n=1}^\infty G_n \zeta^n + G_0 + \sum_{n=1}^\infty K_n \left ({\rho \over \zeta} \right )^n, \\
\text{with} &\qquad a= {d_0 - e_0 \over 4 \log \rho}, \qquad \chi={\rm i} \mu \zeta_{0} [2\text{Im}[z_{0}]-1], \qquad \nu=2{\rm i} \mu \zeta_{0}^{2} \text{Im}[z_{0}].
\end{split}\label{FG SP stresslet}
\end{equation}

\noindent
%Constants $\chi$ and $\nu$ were found by requiring $G(\zeta)$, firstly, to have the correct double pole at $z_{0}$ accounting for the dipole that is associated with the stresslet and, secondly, not having a simple pole at $z_{0}$.
Coefficients $\{F_{n}, H_{n}, G_{n}, K_{n}|n \in \mathbb{N}\}$ and $G_{0}$ are given by \eqref{G0}, \eqref{coef1channel} and \eqref{coef2channel} in terms of coefficients $\{ d_{n}, e_{n} | n \in \mathbb{Z} \}$ which are now given by
%\begin{equation}
%\begin{split}
%\sum_{n=-\infty}^\infty d_n \zeta^n &=-\frac{{\rm i} \overline{\mu} \overline{\zeta_{0}}}{1/\zeta-\overline{\zeta_{0}}}-\frac{\chi}{\zeta-\zeta_{0}}-\frac{\nu}{(\zeta-\zeta_{0})^2}, \\
%\sum_{n=-\infty}^\infty e_n \zeta^n &=-\frac{{\rm i} \overline{\mu} \overline{\zeta_{0}}}{\rho/\zeta-\overline{\zeta_{0}}}-\frac{[{\rm i}\mu\zeta_{0} \log\rho^{2}+\chi]}{\rho \zeta-\zeta_{0}}-\frac{[{\rm i}\mu \zeta_{0}^{2}\log\rho^{2}+\nu]}{(\rho \zeta-\zeta_{0})^2}.
%\end{split}
%\end{equation}
%
%
%\noindent
%Using local expansions, we find
\begin{equation}
\begin{split}
d_{n}&=-{\rm i} \overline{\mu} (\overline{\zeta_{0}})^{n}, \quad n \ge 1, \\
d_{-n}&={\rm i} \mu (\zeta_{0})^{n} [1-2n \text{Im}[z_{0}]], \quad n \ge 1, \\
d_{0}&=0
\end{split}
\end{equation}
and
\begin{equation}
\begin{split} 
e_{n}&=-{\rm i} \mu \rho^{n} (\zeta_{0})^{-n} [1 +n (2\text{Im}[z_{0}]+ \log \rho^{2}) ], \quad n \ge 1, \\
e_{-n}&={\rm i} \overline{\mu} \rho^{n} (\overline{\zeta_{0}})^{-n}, \quad n \ge 1, \\
e_{0}&=2~ \text{Im}[\mu].
\end{split}\label{stressletencoef}
\end{equation}

It turns out that the full analysis just summarized can be bypassed by noticing that the Goursat functions $f(z)$ and $g'(z)$ for the periodic array of stresslets can be derived from those of Stokeslets by computing the parametric derivatives (cf. \cite{DavisCrowdy}):
\begin{equation}\label{param. deriv stresslet}
f(z) \mapsto -\frac{\partial f(z)}{\partial z_{0}}, \qquad g'(z) \mapsto -\frac{\partial g'(z)}{\partial z_{0}}.
\end{equation}
To verify that, we apply the parametric derivatives \eqref{param. deriv stresslet} in \eqref{eq:eqStokeslet} and obtain
\begin{equation}
\begin{split}
f(z)&=\frac{\mu}{z-z_{0}}+ {\rm analytic~function}, \\
g'(z)&=\frac{\mu \overline{z_{0}}}{(z-z_{0})^{2}}-\frac{\overline{\mu}}{z-z_{0}}+ {\rm analytic~function}.
\end{split}\label{local parametric}
\end{equation}
The first terms in $f(z)$ and $g'(z)$ correspond to a stresslet of strength $\mu$ at $z_{0}$ given by \eqref{eq:eqstresslet}, while the second term in $g'(z)$ corresponds to an additional source/sink and/or rotlet singularity. Since our aim is to find expressions for a periodic array of torque-free stresslets without any other induced flow, we can first use \eqref{param. deriv stresslet} in \eqref{FG SP Stokeslet}, where
\begin{equation}
\frac{\partial}{\partial z_{0}}=\frac{\partial \zeta_{0}}{\partial z_{0}} \frac{\partial}{\partial \zeta_{0}}={\rm i} \zeta_{0} \frac{\partial}{\partial \zeta_{0}},
\end{equation}
and, then, subtract terms associated to the simple pole in $g'(z)$. This alternative approach gives identical expressions to \eqref{FG SP stresslet}--\eqref{stressletencoef}.

\subsection{Periodic array of force quadrupoles in a channel}

Consider now the problem of 
a periodic array of force quadrupoles in a channel.
In this problem too the full analysis can be bypassed by noticing that the Goursat functions $f (z)$ and $g'(z)$ for this problem can be derived from those associated with an array of 
torque-free stresslets by computing the parametric derivatives:
\begin{equation}\label{param. deriv force quadrupole}
f(z) \mapsto \frac{\partial f(z)}{\partial z_{0}}, \qquad g'(z) \mapsto \frac{\partial g'(z)}{\partial z_{0}}.
\end{equation}
Following similar steps as previously, or using the parametric derivatives \eqref{param. deriv force quadrupole}, we find that the Goursat functions for a periodic array of force quadrupoles in a channel with local behaviour \eqref{eq:eqforcequadrupole} at $z_{0}$ in the principal window to are given by
\begin{equation}
\begin{split}
F(\zeta) &= \frac{\beta}{\zeta-\zeta_0} +\frac{\gamma}{(\zeta-\zeta_0)^2} + \sum_{n=1}^\infty F_n \zeta^n + \sum_{n=1}^\infty H_n \left ({\rho \over \zeta} \right )^n, \\
G(\zeta) &= \frac{\delta}{\zeta-\zeta_0} + \frac{\epsilon}{(\zeta-\zeta_0)^{2}} + \frac{\kappa}{(\zeta-\zeta_0)^{3}}
+ {\rm i} \log \zeta \left [ {F'(\zeta) \over {\cal Z}'(\zeta)} \right ] \\
&~~~+ \sum_{n=1}^\infty G_n \zeta^n + G_0 + \sum_{n=1}^\infty K_n \left ({\rho \over \zeta} \right )^n, \\
\text{with} & \qquad \beta=-\mu \zeta_{0}, \qquad \gamma=-\mu \zeta_{0}^2, \qquad \delta=2 \mu \zeta_{0} [1-\text{Im}[z_{0}]], \\
&~~~~~~\epsilon=2\mu \zeta_{0}^{2}[1-3\text{Im}[z_{0}]], \qquad \kappa=-4\mu \zeta_{0}^{3} \text{Im}[z_{0}].
\end{split}
\end{equation}
%Constants $\beta, \gamma, \delta, \epsilon$ and $\kappa$ were found by requiring $F(\zeta)$ to have the correct double pole at $z_{0}$ and not having a simple pole at $z_{0}$ and, also, $G(\zeta)$ to have the correct triple pole at $z_{0}$ and not having simple and double poles at $z_{0}$.
The coefficients $\{F_{n}, H_{n}, G_{n}, K_{n}|n \in \mathbb{N}\}$ and $G_{0}$ are again
given by \eqref{G0}, \eqref{coef1channel} and \eqref{coef2channel} in terms of coefficients $\{ d_{n}, e_{n} | n \in \mathbb{Z} \}$ which, in this case, are
%can be found explicitly using local expansions (or, alternatively, numerically using Fast Fourier transforms) in
%\begin{equation}
%\begin{split}
%\sum_{n=-\infty}^\infty d_n \zeta^n &=\frac{\overline{\beta}}{1/\zeta-\overline{\zeta_{0}}}+\frac{\overline{\gamma}}{(1/\zeta-\overline{\zeta_{0}})^2}-\frac{\delta}{\zeta-\zeta_{0}}-\frac{\epsilon}{(\zeta-\zeta_{0})^2}-\frac{\kappa}{(\zeta-\zeta_{0})^3}, \\
%\sum_{n=-\infty}^\infty e_n \zeta^n &=\frac{\overline{\beta}}{\rho/\zeta-\overline{\zeta_{0}}}+\frac{\overline{\gamma}}{(\rho/\zeta-\overline{\zeta_{0}})^2}-\frac{[\beta \log\rho^{2}+\delta]}{\rho \zeta-\zeta_{0}}-\frac{[(\beta \zeta_{0}+2\gamma) \log\rho^{2}+\epsilon]}{(\rho \zeta-\zeta_{0})^2} \\
%&~~~~~-\frac{[2\gamma \zeta_{0} \log\rho^{2}+\kappa]}{(\rho \zeta-\zeta_{0})^3}.
%\end{split}
%\end{equation}
%
%\noindent
%Using local expansions, we find
\begin{equation}
\begin{split}
d_{n}&=-n \overline{\mu} (\overline{\zeta_{0}})^{n}, \quad n \ge 1, \\
d_{-n}&=2n\mu (\zeta_{0})^{n} [n \text{Im}[z_{0}]-1], \quad n \ge 1, \\
d_{0}&=0
\end{split}
\end{equation}
and
\begin{equation}
\begin{split} 
e_{n}&=-n \mu \rho^{n} (\zeta_{0})^{-n} [2+n( 2\text{Im}[z_{0}]+\log \rho^{2})], \quad n \ge 1, \\
e_{-n}&=-n \overline{\mu} \rho^{n} (\overline{\zeta_{0}})^{-n}, \quad n \ge 1, \\
e_{0}&=0.
\end{split}
\end{equation}

\subsection{Periodic array of dipoles and quadrupoles in a channel \label{Sec8}}

The Goursat functions $f (z)$ and $g'(z)$ for a periodic array of dipoles can be derived from those of Stokeslets \eqref{FG SP Stokeslet} by computing the mixed parametric derivatives:
\begin{equation}\label{mixed param. deriv dipole}
f(z) \mapsto \frac{\partial^{2} f(z)}{\partial z_{0} \partial \overline{z_{0}}}, \qquad g'(z) \mapsto \frac{\partial^{2} g'(z)}{\partial z_{0} \partial \overline{z_{0}}}.
\end{equation}
where
\begin{equation}
\frac{\partial^{2}}{\partial z_{0} \partial \overline{z_{0}}}=\zeta_{0} \overline{\zeta_{0}} \frac{\partial^{2}}{\partial \zeta_{0} \partial \overline{\zeta_{0}}}.
\end{equation}
The Goursat functions for a periodic array of dipoles in a channel with local behaviour \eqref{sourcedipole} at $z_{0}$ in the principal period window are
\begin{equation}
\begin{split}
F(\zeta) &= \sum_{n=1}^\infty F_n \zeta^n + \sum_{n=1}^\infty H_n \left ({\rho \over \zeta} \right )^n, \\
G(\zeta) &= \frac{\mu \zeta_{0}}{\zeta-\zeta_0} + \frac{\mu \zeta_{0}^2}{(\zeta-\zeta_0)^{2}}
+ {\rm i} \log \zeta \left [ {F'(\zeta) \over {\cal Z}'(\zeta)} \right ] + \sum_{n=1}^\infty G_n \zeta^n + G_0 + \sum_{n=1}^\infty K_n \left ({\rho \over \zeta} \right )^n.
\end{split}
\end{equation}

\noindent
%Constants $\chi$ and $\lambda$ were found by requiring $G(\zeta)$ to have the correct double pole at $z_{0}$ and not having a simple pole at $z_{0}$.
Again, the coefficients $\{F_{n}, H_{n}, G_{n}, K_{n}|n \in \mathbb{N}\}$ and $G_{0}$ are given by \eqref{G0}, \eqref{coef1channel} and \eqref{coef2channel} in terms of coefficients $\{ d_{n}, e_{n} | n \in \mathbb{Z} \}$ which are given by
\begin{equation}
\begin{split}
d_{n}&=0, \quad n \ge 1, \\
d_{-n}&=-n \mu (\zeta_{0})^{n}, \quad n \ge 1, \\
d_{0}&=0
\end{split}
\end{equation}
and
\begin{equation}
\begin{split} 
e_{n}&=-n \mu \rho^{n} (\zeta_{0})^{-n}, \quad n \ge 1, \\
e_{-n}&=0, \quad n \ge 1, \\
e_{0}&=0.
\end{split}
\end{equation}

Similarly the Goursat functions $f (z)$ and $g'(z)$ for a periodic array of quadrupoles can be derived from those of torque-free stresslets \eqref{FG SP stresslet} by computing the mixed parametric derivatives:
\begin{equation}\label{mixed param. deriv quadrupole}
f(z) \mapsto -\frac{\partial^{2} f(z)}{\partial z_{0} \partial \overline{z_{0}}}, \qquad g'(z) \mapsto -\frac{\partial^{2} g'(z)}{\partial z_{0} \partial \overline{z_{0}}}.
\end{equation}

\noindent
The Goursat functions for a periodic array of quadrupoles in a channel with local behaviour \eqref{sourcequadrupole} at $z_{0}$ in the principal period window are
\begin{equation}
\begin{split}
F(\zeta) &= \sum_{n=1}^\infty F_n \zeta^n + \sum_{n=1}^\infty H_n \left ({\rho \over \zeta} \right )^n, \\
G(\zeta) &= \frac{{\rm i}\mu \zeta_{0}}{\zeta-\zeta_0} + \frac{3{\rm i} \mu \zeta_{0}^{2}}{(\zeta-\zeta_0)^{2}} + \frac{2{\rm i}\mu \zeta_{0}^{3}}{(\zeta-\zeta_0)^{3}}
+ {\rm i} \log \zeta \left [ {F'(\zeta) \over {\cal Z}'(\zeta)} \right ] \\
&~~~+ \sum_{n=1}^\infty G_n \zeta^n + G_0 + \sum_{n=1}^\infty K_n \left ({\rho \over \zeta} \right )^n
\end{split}
\end{equation}

\noindent
%Constants $\delta, \epsilon$ and $\kappa$ were found by requiring $G(\zeta)$ to have the correct triple pole at $z_{0}$ and not having simple and double poles at $z_{0}$.
with coefficients $\{F_{n}, H_{n}, G_{n}, K_{n}|n \in \mathbb{N}\}$ and $G_{0}$ given by \eqref{G0}, \eqref{coef1channel} and \eqref{coef2channel} in terms of coefficients $\{ d_{n}, e_{n} | n \in \mathbb{Z} \}$:
\begin{equation}
\begin{split}
d_{n}&=0, \quad n \ge 1, \\
d_{-n}&=-{\rm i} n^{2} \mu (\zeta_{0})^{n}, \quad n \ge 1, \\
d_{0}&=0
\end{split}
\end{equation}
and
\begin{equation}
\begin{split} 
e_{n}&={\rm i}n^{2} \mu \rho^{n} (\zeta_{0})^{-n}, \quad n \ge 1, \\
e_{-n}&=0, \quad n \ge 1, \\
e_{0}&=0.
\end{split}
\end{equation}

\section{Periodic singularity arrays in the half-plane}

In the limit where the upper wall of the channel becomes distant from the lower wall
we encounter the situation of a periodic array of Stokes singularities in the upper half
plane near an infinite
straight wall and it is useful to consider this situation too since it arises in many physical
modelling situations (e.g. \cite{MannanCortez}).

\begin{figure}
\begin{center}
\includegraphics[scale=1]{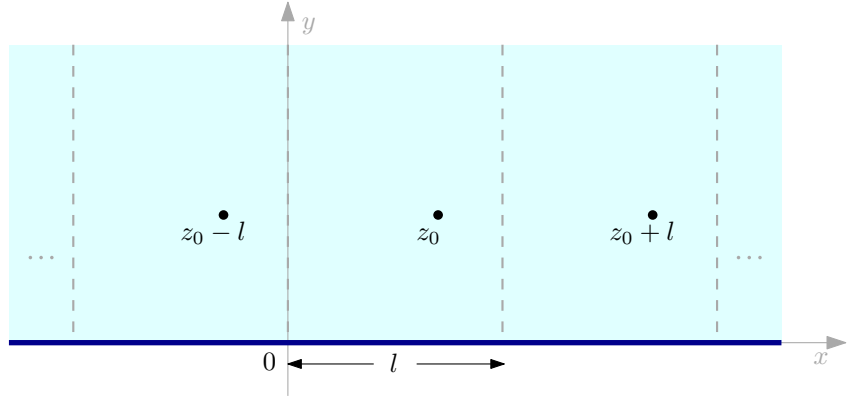}
\caption{A periodic array of point singularities at points  $z=z_{0}+n l$, $n \in \mathbb{Z}$ in the upper half-plane $-\infty<x<\infty$, $y \geq 0$.}
\label{fig: confupper}
\end{center}
\end{figure}

Consider a periodic array of point singularities at $z = z_0 + nl$, $n \in \mathbb{Z}$, with $0 < \text{Re}[z_0] < l$, $\text{Im}[z_0] > 0$, in the upper half-plane $-\infty<x<\infty$, $y \geq 0$. Figure \ref{fig: confupper} shows a schematic of the configuration. While the problem can be solved using standard Fourier transform techniques \cite{Poz},  we give an alternative derivation, and form, of the solution by extending the ideas given in the previous sections.
%In this section we report the Goursat function representations for singly periodic singularity arrays in the upper half-plane (Fig. \ref{fig: confupper}); these follow from those presented in \S \ref{Sec5} by taking the limit $h \rightarrow \infty$ which corresponds to $\rho \rightarrow 0$.

\subsection{Periodic array of Stokeslets in the upper half-plane}

Consider a periodic array of Stokeslets each of strength given by \eqref{Strength} for some $\mu \in \mathbb{C}$. The representative Stokeslet in the fundamental period semi-strip shown in Fig. \ref{fig: confupper} is located at $z_0$ and has a preimage at $\zeta_0$ with $z_0=z(\zeta_0)=-{\rm i}\log \zeta_0$. 
The same conformal mapping (\ref{eq:eqmap}) is employed, but now we take the limit
$\rho \to 0$ so that the preimage domain in the unit disc in the parametric $\zeta$ plane.
The unit circle $|\zeta|=1$ is transplanted to the no-slip wall along the $x$ axis.

Again we introduce the functions $F(\zeta)$ and $G(\zeta)$ defined in (\ref{FGzeta}).
%\begin{equation}
%F(\zeta) \equiv f({\cal Z}(\zeta)), \qquad G(\zeta) \equiv g'({\cal Z}(\zeta)).
%\end{equation}
Let
\begin{equation}
F(\zeta)=\mu \log(\zeta-\zeta_0)+\epsilon \log(1- \overline{\zeta_0} \zeta)+\frac{\kappa}{1- \overline{\zeta_0} \zeta}+\lambda,
\end{equation}
where constants $\epsilon$, $\kappa$ and $\lambda$ are to be found. 
In terms of functions of the $\zeta$ variable
the no-slip boundary condition on the wall $\overline{z}=z$ takes the form
\begin{equation}
-\overline{F}(1/\zeta)+{\cal Z}(\zeta) \left[ \frac{F'(\zeta)}{{\cal Z}'(\zeta)} \right]+G(\zeta)=0,
\end{equation}
where we have that $\overline{\zeta}=1/\zeta$ on $|\zeta|=1$; this implies that
\begin{equation}
G(\zeta)=\overline{F}(1/\zeta)-{\cal Z}(\zeta) \left[ \frac{F'(\zeta)}{{\cal Z}'(\zeta)} \right].
\end{equation}
It remains to find the unknown constants $\epsilon$, $\kappa$ and $\lambda$. These follow by insisting that $u-{\rm i}v$ is invariant as $\zeta\mapsto \zeta \mathrm{e}^{2\pi{\rm i}}$ which corresponds to the periodicity of the velocity field and by ensuring that the singularity of $u-{\rm i}v$ at $\zeta_0$ has the required form. It is found that
\begin{equation}
\epsilon=-\mu, \qquad \kappa=\overline{\mu} \log|\zeta_0|^2, \qquad \lambda=-\kappa.
\end{equation}

In summary, the Goursat functions for a periodic array of Stokeslets in the upper half-plane with one at $z_{0}$ in the principal period window and having strength given by \eqref{Strength} are
\begin{equation}
\begin{split}
F(\zeta) &= \mu \log(\zeta-\zeta_0)-\mu \log(1- \overline{\zeta_0} \zeta)+\frac{\kappa}{1- \overline{\zeta_0} \zeta}-\kappa, \\
G(\zeta) &= - \overline{\mu} \log (\zeta-\zeta_0) + \overline{\mu} \log(1- \overline{\zeta_0} \zeta) + {\overline{\kappa} \zeta_0 \over {\zeta-\zeta_0}} + {\rm i} \log \zeta \left [ {F'(\zeta) \over {\cal Z}'(\zeta)} \right ].
\end{split}\label{Stokesletsupper}
\end{equation}
As expected, these expressions 
coincide with those obtained by taking $\rho \to 0$ 
in the solutions presented in \S \ref{Sec5} (which corresponds to the limit $h \rightarrow \infty$).

%\begin{equation}
%\begin{split}
%F(\zeta) &= \mu \log (\zeta-\zeta_0) + \sum_{n=1}^\infty F_n \zeta^n, \\
%G(\zeta) &= - \overline{\mu} \log (\zeta-\zeta_0) + {\lambda \over \zeta-\zeta_0}
%+ {\rm i} \log \zeta \left [ {F'(\zeta) \over z'(\zeta)} \right ] + \sum_{n=1}^\infty G_n \zeta^n, \\
%\text{with} &\qquad \lambda = \mu \zeta_0 \log |\zeta_0|^2,
%\end{split}\label{Stokesletsupper}
%\end{equation}
%where coefficients $\{F_{n}, G_{n}|n \in \mathbb{N}\}$ are given by
%\begin{equation}
%F_{n} = (\overline{\zeta_{0}})^{n}(\mu/n+\overline{\mu} \log|\zeta_{0}|^{2}), \qquad G_n = -\overline{\mu} (\overline{\zeta_{0}})^{n}/n, \qquad n \geq 1.
%\end{equation}

\subsection{Periodic array of stresslets in the upper half-plane}

Following similar steps or taking parametric derivatives as shown in \S \ref{Sec5}, we find that the Goursat functions for a periodic array of torque-free stresslets in the upper half-plane with local behaviour \eqref{eq:eqstresslet} at $z_{0}$ in the principal period window are
\begin{equation}
\begin{split}
F(\zeta) &= \frac{{\rm i} \mu \zeta_{0}}{\zeta-\zeta_0} + \frac{\beta}{1-\overline{\zeta_0} \zeta} + \frac{\gamma}{(1-\overline{\zeta_0} \zeta)^2}-{\rm i}\overline{\mu}, \\
G(\zeta) &= \frac{\chi}{\zeta-\zeta_0} + \frac{\nu}{(\zeta-\zeta_0)^{2}} -\frac{{\rm i}\overline{\mu}}{1-\overline{\zeta_0}\zeta}
+ {\rm i} \log \zeta \left [ {F'(\zeta) \over {\cal Z}'(\zeta)} \right ] + {\rm i}\overline{\mu}, \\
\text{with}&~~~~~~\beta={\rm i}\overline{\mu}[2\text{Im}[z_0]+1], \qquad \gamma=-2{\rm i}\overline{\mu}\text{Im}[z_0], \\
&~~~~~~\chi={\rm i} \mu \zeta_{0} [2\text{Im}[z_{0}]-1], \qquad \nu=2{\rm i} \mu \zeta_{0}^{2} \text{Im}[z_{0}].
\end{split}
\end{equation}
%\begin{equation}
%\begin{split}
%F(\zeta) &= \frac{{\rm i} \mu \zeta_{0}}{\zeta-\zeta_0} + \sum_{n=1}^\infty F_n \zeta^n, \\
%G(\zeta) &= \frac{\chi}{\zeta-\zeta_0} + \frac{\nu}{(\zeta-\zeta_0)^{2}}
%+ {\rm i} \log \zeta \left [ {F'(\zeta) \over z'(\zeta)} \right ] + \sum_{n=1}^\infty G_n \zeta^n, \\
%\text{with}&~~~~~~\chi={\rm i} \mu \zeta_{0} [2\text{Im}[z_{0}]-1], \qquad \nu=2{\rm i} \mu \zeta_{0}^{2} \text{Im}[z_{0}],
%\end{split}
%\end{equation}
%where
%\begin{equation}
%F_{n} = {\rm i} \overline{\mu} (\overline{\zeta_{0}})^{n} [1-2n \text{Im}[z_{0}]], \qquad G_n = -{\rm i} \overline{\mu} (\overline{\zeta_{0}})^{n}, \qquad n \geq 1. 
%\end{equation}

\subsection{Periodic array of force quadrupoles in the upper half-plane}

The Goursat functions for a periodic array of force quadrupoles in the upper half-plane with local behaviour \eqref{eq:eqforcequadrupole} at $z_{0}$ in the principal period window are
\begin{equation}
\begin{split}
F(\zeta) &= -\frac{\mu \zeta_0}{\zeta-\zeta_0} -\frac{\mu \zeta_0^2}{(\zeta-\zeta_0)^2} +\frac{\beta}{1-\overline{\zeta_0}\zeta}+\frac{\gamma}{(1-\overline{\zeta_0}\zeta)^2}+\frac{\delta}{(1-\overline{\zeta_0}\zeta)^3}, \\
G(\zeta) &= \frac{\epsilon}{\zeta-\zeta_0} + \frac{\kappa}{(\zeta-\zeta_0)^{2}} + \frac{\lambda}{(\zeta-\zeta_0)^{3}} +\frac{\overline{\mu}}{1-\overline{\zeta_0} \zeta}-\frac{\overline{\mu}}{(1-\overline{\zeta_0} \zeta)^2} + {\rm i} \log \zeta \left [ {F'(\zeta) \over {\cal Z}'(\zeta)} \right ], \\
\text{with} & \qquad \beta=-2\overline{\mu}[1+\text{Im}[z_0]], \qquad \gamma=2\overline{\mu}[1+3\text{Im}[z_0]], \qquad \delta=-4\overline{\mu} \text{Im}[z_0], \\
&\qquad \epsilon=2 \mu \zeta_{0} [1-\text{Im}[z_{0}]], \qquad \kappa=2\mu \zeta_0^2[1-3\text{Im}[z_0]], \qquad \lambda=-4\mu \zeta_0^3 \text{Im}[z_0].
\end{split}
\end{equation}
%\begin{equation}
%\begin{split}
%F(\zeta) &= \frac{\beta}{\zeta-\zeta_0} +\frac{\gamma}{(\zeta-\zeta_0)^2} + \sum_{n=1}^\infty F_n \zeta^n, \\
%G(\zeta) &= \frac{\delta}{\zeta-\zeta_0} + \frac{\epsilon}{(\zeta-\zeta_0)^{2}} + \frac{\kappa}{(\zeta-\zeta_0)^{3}}
%+ {\rm i} \log \zeta \left [ {F'(\zeta) \over z'(\zeta)} \right ] + \sum_{n=1}^\infty G_n \zeta^n, \\
%\text{with} & \qquad \beta=-\mu \zeta_{0}, \qquad \gamma=-\mu \zeta_{0}^2, \qquad \delta=2 \mu \zeta_{0} [1-\text{Im}[z_{0}]], \\
%&~~~~~~\epsilon=2\mu \zeta_{0}^{2}[1-3\text{Im}[z_{0}]], \qquad \kappa=-4\mu \zeta_{0}^{3} \text{Im}[z_{0}],
%\end{split}
%\end{equation}
%where
%\begin{equation}
%F_{n} =-2n \overline{\mu} (\overline{\zeta_{0}})^{n} [n \text{Im}[z_{0}]-1],  \qquad G_n = -n \overline{\mu} (\overline{\zeta_{0}})^{n}, \qquad n \geq 1.
%\end{equation}

\subsection{Periodic array of dipoles and quadrupoles in the upper half-plane \label{Sec8}}

The Goursat functions for a periodic array of dipoles in the upper half-plane with local behaviour \eqref{sourcedipole} at $z_{0}$ in the principal period window are
\begin{equation}
\begin{split}
F(\zeta) &= -\frac{\overline{\mu}}{1-\overline{\zeta_0}\zeta}+\frac{\overline{\mu}}{(1-\overline{\zeta_0}\zeta)^2}, \\
G(\zeta) &= \frac{\mu \zeta_{0}}{\zeta-\zeta_0} + \frac{\mu \zeta_{0}^2}{(\zeta-\zeta_0)^{2}}
+ {\rm i} \log \zeta \left [ {F'(\zeta) \over {\cal Z}'(\zeta)} \right ].
\end{split}
\end{equation}
%\begin{equation}
%\begin{split}
%F(\zeta) &= \sum_{n=1}^\infty F_n \zeta^n, \\
%G(\zeta) &= \frac{\mu \zeta_{0}}{\zeta-\zeta_0} + \frac{\mu \zeta_{0}^2}{(\zeta-\zeta_0)^{2}}
%+ {\rm i} \log \zeta \left [ {F'(\zeta) \over z'(\zeta)} \right ],
%\end{split}
%\end{equation}
%where
%\begin{equation}
%F_{n} =n \overline{\mu} (\overline{\zeta_{0}})^{n}, \quad n \geq 1.
%\end{equation}

The Goursat functions for a periodic array of quadrupoles in the upper half-plane with local behaviour \eqref{sourcequadrupole} at $z_{0}$ in the principal period window are
\begin{equation}
\begin{split}
F(\zeta) &= -\frac{{\rm i}\overline{\mu}}{1-\overline{\zeta_0}\zeta}+\frac{3{\rm i}\overline{\mu}}{(1-\overline{\zeta_0}\zeta)^2}-\frac{2{\rm i}\overline{\mu}}{(1-\overline{\zeta_0}\zeta)^3}, \\
G(\zeta) &= \frac{{\rm i}\mu \zeta_{0}}{\zeta-\zeta_0} + \frac{3{\rm i} \mu \zeta_{0}^{2}}{(\zeta-\zeta_0)^{2}} + \frac{2{\rm i}\mu \zeta_{0}^{3}}{(\zeta-\zeta_0)^{3}}
+ {\rm i} \log \zeta \left [ {F'(\zeta) \over {\cal Z}'(\zeta)} \right ].
\end{split}
\end{equation}

\section{Discussion}

New analytical representations for a wide range of Stokes flows due to periodic arrays of point singularities in a two-dimensional no-slip channel and in a half-plane near a wall have been given. 
The solutions are given explicitly as functions of a parametric $\zeta$ variable.
If preferred,
in all cases, the solutions can be expressed directly in terms of the variable
$z=x+{\rm i}y$ on use of (\ref{zetaz}).

We have shown that the associated Laurent expansion terms (for the channel geometry) decay rapidly which suggests that only a few terms are, in general, required to provide accurate and efficient computations. The analytical expressions presented
here can be used to model Stokes flow problems in singly periodic geometries; for example, one can model a flow past a periodic array of wall-attached straight barriers in a channel or in the half-plane \cite{Davis1993}, or model cilia sheets \cite{MannanCortez}.

\vspace{0.3cm}
\noindent
{\bf Acknowledgements} ~DGC was supported by an EPSRC Established Career Fellowship (EP/K019430/10) and by a Royal Society Wolfson Research Merit Award. Both authors acknowledge financial support from a Research Grant from the Leverhulme Trust.

\appendix

\section{Transform method \label{appendixA}}

We briefly discuss how to extend a novel transform approach to biharmonic boundary value problems for both polygonal and circular domains recently described by the authors in \cite{LucaCrowdy} to reappraise the problem of a periodic array of point singularities in a channel geometry thereby offering an alternative to the solutions given in \S \ref{Sec5}.

The Goursat functions can be represented by
\begin{equation}\label{expansions}
f(z)=f_{s}(z)+f_{R}(z), \qquad g'(z)=g'_{s}(z)+g'_{R}(z),
\end{equation}
where $f_{s}(z)$, $g'_{s}(z)$ are related to the point singularity at $z_{0}$ (for example, for a point stresslet at $z_{0}$, they are given by \eqref{f stresslet} and \eqref{g' stresslet}) and $f_{R}(z)$, $g'_{R}(z)$ are the correction functions to be found. The functions $f_{R}(z)$, $g'_{R}(z)$ are analytic and single-valued in the fluid region and they have the following integral representations \cite{LucaCrowdy}:
\begin{equation}\label{representations1}
f_{R}(z)=\frac{1}{2 \pi} \left[ \sum_{j=1}^{4} \int_{L_j} { \rho_{j}(k) \mathrm{e}^{{\rm i} k z} \mathrm{d}k} \right], \qquad g'_{R}(z)=\frac{1}{2 \pi} \left[ \sum_{j=1}^{4} \int_{L_j} { \hat{\rho}_{j}(k) \mathrm{e}^{{\rm i} k z} \mathrm{d}k} \right],
\end{equation}
where $L_j$, $j=1,2,3,4$ are oriented rays from 0 in the spectral $k$-plane \cite{LucaCrowdy} and $\rho_{j}(k), \hat{\rho}_{j}(k)$, $j=1,2,3,4$ are the spectral functions defined by
\begin{equation}\label{representations2}
\begin{split}
\rho_{1}(k)&=\int_{0}^{l} {f_{R}(z) \mathrm{e}^{-{\rm i} k z} \mathrm{d}z}, \qquad ~~~\rho_{2}(k)=\int_{l}^{l+{\rm i}h} {f_{R}(z) \mathrm{e}^{-{\rm i} k z} \mathrm{d}z}, \\
\rho_{3}(k)&=\int_{l+{\rm i}h}^{{\rm i}h} {f_{R}(z) \mathrm{e}^{-{\rm i} k z} \mathrm{d}z}, \qquad \rho_{4}(k)=\int_{{\rm i}h}^{0} {f_{R}(z) \mathrm{e}^{-{\rm i} k z} \mathrm{d}z},
\end{split}
\end{equation}
and
\begin{equation}\label{representations2g}
\begin{split}
\hat{\rho}_{1}(k)&=\int_{0}^{l} {g'_{R}(z) \mathrm{e}^{-{\rm i} k z} \mathrm{d}z}, \qquad ~~~ \hat{\rho}_{2}(k)=\int_{l}^{l+{\rm i}h} {g'_{R}(z) \mathrm{e}^{-{\rm i} k z} \mathrm{d}z}, \\
\hat{\rho}_{3}(k)&=\int_{l+{\rm i}h}^{{\rm i}h} {g'_{R}(z) \mathrm{e}^{-{\rm i} k z} \mathrm{d}z}, \qquad \hat{\rho}_{4}(k)=\int_{{\rm i}h}^{0} {g'_{R}(z) \mathrm{e}^{-{\rm i} k z} \mathrm{d}z}.
\end{split}
\end{equation}
The spectral functions satisfy the so-called global relations:
\begin{equation}\label{GR}
\sum_{j=1}^{4} {\rho_{j}(k)}=0, \qquad \sum_{j=1}^{4} {\hat{\rho}_{j}(k)}=0, \qquad \text{for } k \in \mathbb{C}.
\end{equation}

The analysis of the boundary and periodicity conditions allows us to deduce relations between the spectral functions. We omit the details and report the key expressions; these are:

\begin{align}
-\overline{\rho_{1}}(-k) -\frac{\partial [k \rho_{1}(k)]}{\partial k}+\hat{\rho}_{1}(k)+ l f_{R}(l) \mathrm{e}^{-{\rm i}kl}&=R_{1}(k), \label{finalTM1} \\
\rho_{4}(k)+\mathrm{e}^{{\rm i}kl} \rho_{2}(k)+d ~ q(k)&=R_{2}(k), \label{finalTM2} \\
-\mathrm{e}^{2kh} \overline{\rho_{3}}(-k)-\frac{\partial [k \rho_{3}(k)]}{\partial k}+2kh \rho_{3}(k)+\hat{\rho}_{3}(k) +r(k)&=R_{3}(k), \label{finalTM3} \\
-{\rm i}k l \rho_{4}(k)+\hat{\rho}_{4}(k)+\mathrm{e}^{{\rm i}k l} \hat{\rho}_{2}(k)-l f_{R}(0)+l f_{R}({\rm i}h) \mathrm{e}^{kh} +\overline{d} ~ q(k)&=R_{4}(k), \label{finalTM4}
\end{align}
where $d \in \mathbb{C}$ is a constant and
\begin{align}
R_{1}(k) &\equiv \int_{0}^{l} {[\overline{f_{s}(z)}-z f'_{s}(z)-g'_{s}(z)] \mathrm{e}^{-{\rm i}kz} \mathrm{d}z}, \\
R_{2}(k) &\equiv \int_{{\rm i}h}^{0} {[-f_{s}(z)+f_{s}(z+l)] \mathrm{e}^{-{\rm i}kz} \mathrm{d}z}, \\
R_{3}(k) &\equiv \int_{l+{\rm i}h}^{{\rm i}h} {[\overline{f_{s}(z)}-(z-2{\rm i}h) f'_{s}(z)-g'_{s}(z)] \mathrm{e}^{-{\rm i}kz} \mathrm{d}z}, \\
R_{4}(k) &\equiv \int_{{\rm i}h}^{0} {[l f'_{s}(z) - g'_{s}(z)+g'_{s}(z+l)] \mathrm{e}^{-{\rm i}kz} \mathrm{d}z}
\end{align}
and
\begin{equation}
q(k) \equiv \int_{{\rm i}h}^{0} {\mathrm{e}^{-{\rm i} kz} \mathrm{d}z}, \qquad r(k)=-{\rm i}h f_{R}({\rm i}h) \mathrm{e}^{kh}-(l-{\rm i}h) f_{R}(l+{\rm i}h) \mathrm{e}^{-{\rm i}k(l+{\rm i}h)}.
\end{equation}
\\

\noindent
Addition of \eqref{finalTM1} and \eqref{finalTM3} and use of \eqref{GR},\eqref{finalTM2},\eqref{finalTM4} gives, after some algebra,
\begin{equation}\label{rho1 expression}
\rho_{1}(k)=\frac{2kh W(k)-(\mathrm{e}^{2kh}-1) \overline{W}(-k)}{4 [\sinh^2(kh)-k^2 h^2]}, 
\end{equation}
where $W(k)$ contains $\rho_{4}(k), \hat{\rho}_{4}(k), f_R(0)$, $f_R({\rm i}h), d$ and known quantities. The spectral function $\rho_{1}(k)$ is analytic everywhere in the complex $k$-plane which means that its numerator in \eqref{rho1 expression} must vanish at zeros of its denominator in the $k$-plane, i.e. we must require
\begin{equation} \label{conditions 1}
2kh W(k)-(\mathrm{e}^{2kh}-1) \overline{W}(-k)=0, \qquad \text{for} \quad k \in \Sigma_{1} \equiv \{k \in \mathbb{C}| \sinh^2(kh)-k^2 h^2]=0\},
\end{equation}
together with conditions at $k=0$ following from \eqref{rho1 expression}. Next, we use the series representations
\begin{equation}\label{sums expansions}
f_{R}(z)= \sum_{m} {a_{m} T_m(z)}, \quad g_{R}(z)=\sum_{m} {b_{m} T_m(z)}, \quad \text{along} \quad z={\rm i}y, y\in[0,h],
\end{equation}
where $a_{m}, b_{m} \in \mathbb{C}$ are unknown coefficients and $T_m(z)$ are basis functions (e.g. Fourier, Chebyshev), and truncate the sums in \eqref{sums expansions} to finite number of terms. Then, we formulate a linear system for the unknown coefficients $a_{m}$, $b_{m}$, parameter $d$ and their complex conjugates. The linear system comprises conditions \eqref{conditions 1} evaluated at as many points in the set $\Sigma_{1}$ as needed, together with conditions at $k=0$. Once the unknowns are computed from the solution of the truncated linear system, the spectral functions $\rho_{4}(k)$ and $\hat{\rho}_{4}(k)$ can be found. The remaining spectral functions can be found by back substitution into \eqref{finalTM1}--\eqref{finalTM4}, and therefore the correction functions $f_{R}(z)$ and $g'_{R}(z)$ can be computed.

\end{document}